\documentclass{amsart}
\usepackage[dvips]{graphicx}

\newtheorem{theorem}{Theorem}[section]

\newtheorem{proposition}[theorem]{Proposition}
\newtheorem{corollary}[theorem]{Corollary}
\theoremstyle{definition}
\newtheorem{definition}[theorem]{Definition}
\newtheorem{remark}[theorem]{Remark}

\numberwithin{equation}{section}

\newcommand{\cpk}{\mathbb{P}^4}
\newcommand{\zz}{\mathbb{Z}}
\newcommand{\cc}{\mathbb{C}}
\newcommand{\rr}{\mathbb{R}}
\newcommand{\qq}{\mathbb{Q}}

\newcommand{\hh}{\mathbb{H}}

\newcommand{\sss}{\mathbb{S}}
\newcommand{\lan}{\langle}
\newcommand{\ran}{\rangle}

\newcommand{\Ns}{N_{\sigma}}
\newcommand{\Gs}{G_{\sigma}}
\newcommand{\Gt}{G_{\tau}}

\newcommand{\Us}{U_{\sigma}}
\newcommand{\Usp}{U_{\sigma'}}
\newcommand{\Ut}{U_{\tau}}
\newcommand{\Ms}{M_{\sigma}}

\newcommand{\csigp}{\check{\sigma}'}

\newcommand{\kG}{K({\bf g})}

\begin{document}

\title[The Chen-Ruan Cohomology Ring of 
Mirror Quintic]{The Chen-Ruan Cohomology Ring of \\
Mirror Quintic}

\author{B. Doug Park}
\address{Department of Pure Mathematics, University of Waterloo,
Waterloo, Ontario N2L 3G1, Canada}
\email{bdpark@math.uwaterloo.ca}

\author{Mainak Poddar}
\address{Department of Mathematics,
Michigan State University, E. Lansing, MI 48824, USA}
\email{poddar@math.msu.edu}


\date{July 4, 2002.  Revised on October 10, 2002.}

\begin{abstract}
We compute the Chen-Ruan orbifold cohomology ring of the
Batyrev mirror orbifold of a smooth quintic hypersurface in $\cpk$.  
We identify the obstruction bundle for this example by 
using the Riemann bilinear relations for periods.
We outline a general method of computing the 
Chen-Ruan ring for Calabi-Yau hypersurfaces in projective 
simplicial toric varieties, modulo a conjecture that the
Riemann bilinear relations are adequate for identifying the
obstruction bundle for any complex orbifold.  
\end{abstract}

\maketitle

\section{Introduction}

The Chen-Ruan orbifold cohomology (cf.$\;$\cite{[CR1]}) gives mathematical meaning,
\`a la Gromov-Witten theory, to the various invariants such as ``orbifold Euler characteristic''
and ``orbifold cohomology'' for global quotient orbifolds that have been in use, courtesy of 
string theorists, for quite a while.  The most interesting feature of this new cohomology 
theory, besides the generalization to non-global quotients, is the existence of a ring 
structure which was previously missing.  The ring structure is obtained from Chen-Ruan's 
orbifold quantum cohomology construction (cf.$\;$\cite{[CR2]}) by restricting to the 
class of maps called ghost maps, in the same spirit as the ordinary cup product may be 
recovered from quantum cup product.  

One of the motivations of Chen and Ruan was to construct a mathematical theory that 
is rich enough to formulate mirror symmetry for Calabi-Yau orbifolds without having to 
resort to their smooth crepant resolutions which do not always exist in dimensions 
higher than three.  Although the Chen-Ruan cohomology has had several successes, it is 
yet to make serious progress in this direction.  One problem has been that the 
computation of the ring structure was an unsolved problem when the ``obstruction bundle''
is not trivial, which is indeed the case for the Batyrev mirror families of hypersurfaces
in Fano toric varieties.  

Here we try to rectify the situation by suggesting that the Riemann bilinear relations 
for periods \cite{[Gr]} may be adequate to identify the obstruction bundle.  Although
we do not prove this, we show that this is true in the mirror 
quintic example.  We then proceed to give an explicit description of the ring structure for  
this example.  The integral of the Euler class of the obstruction bundle is computed 
first by an ad hoc method that is shorter but rather restrictive, and then recomputed by using
localization techniques which should work for any toric Calabi-Yau hypersurface.  
We deliberately avoid making use of the fact that the mirror quintic is a global quotient.  
Rather, we focus on developing a strategy that should work in the general setting.  
For instance the twisted sectors of the mirror quintic can be readily determined since  
it is a global quotient. However we obtained them as a special case of a more general theorem.
We would like to point out that the computation of the ordinary cup product of
ample toric
hypersurfaces is made possible by recent results of Mavlyutov \cite{[Ma]} and 
 Gr\"obner basis method (cf.$\;$\cite{[CLO]}). 

Finally, we draw the reader's attention to a relevant conjecture of Ruan.
String theory suggests that the orbifold quantum cohomology ring of a Calabi-Yau orbifold 
should be isomorphic to the ordinary quantum cohomology ring of its crepant 
resolution. This is, of course, not easy to check. 
It is known that the Chen-Ruan ring structure is not preserved under (partial)
crepant resolutions. But the Cohomological Crepant Resolution Conjecture (cf.$\;$\cite{[Ru2]})
specifies how the ordinary cup product of the crepant resolution 
may be deformed by using certain quantum corrections 
to make it equal to the orbifold cup product of the original space. 
Developing the machinery to compute the Gromov-Witten invariants involved in these
correction terms seems to be an interesting problem.  
The reader may find a parallel approach to stringy cohomolgy spaces in \cite{[BM]}
interesting as well.

\subsection*{Acknowledgments}
We thank Alejandro Adem, Karl Heinz Dovermann, Ronald Fintushel, John McCarthy, 
Yongbin Ruan
and Jon Wolfson for very helpful discussions and encouragements. We especially
thank Yongbin Ruan for explaining to us how localization might be used to compute
integral of the Euler class of the obstruction bundle.
The second author also benefitted from conversations with Yuan-Pin Lee, Ernesto Lupercio 
and Bernardo Uribe.  Some computations were done with the aid of {\sl
Mathematica}$^{\circledR}$   Version 4.1 (cf.$\;$\cite{[Wo]}).

\smallskip

\section{Preliminaries}
\subsection{Mirror Quintic}\label{subsec:quintic}

Let $\zeta = \exp (2\pi i /5)$.  Consider the following action of $(\zz_5)^3$ on
$\cpk$:
$$ [x_1 : x_2 :x_3:x_4:x_5] \mapsto [ \zeta^{a_1} x_1 :  \zeta^{a_2} x_2 :
\zeta^{a_3} x_3: \zeta^{a_4} x_4: \zeta^{a_5} x_5] $$
where $a_i \in \zz_5\,$  and  $\,\sum_{i=1}^{5} a_i \equiv 0 \mod 5$.  
From now on, an element of $(\zz_5)^3$ will be denoted by a 5-tuple\/  
$g = ( \zeta^{a_1}, \zeta^{a_2}, \zeta^{a_3}, \zeta^{a_4}, \zeta^{a_5})$\/ 
satisfying the above congruence.

The mirror family (cf.$\;$\cite{[B1]},$\,$\cite{[GP]})
of smooth quintic hypersurfaces in $\cpk$ is given by a
one-parameter family of quintic hypersurfaces 
$X_{\psi} \subset \cpk/(\zz_5)^3$ satisfying the
following equations:
\begin{equation}\label{eq:quintic}
x_1^5 + x_2^5 + x_3^5 + x_4^5 + x_5^5 + \psi \,x_1 x_2
x_3 x_4 x_5 \:=\: 0 \, .
\end{equation}
Every member of this mirror family is a
3-dimensional Calabi-Yau orbifold, provided that $\psi \neq -5\zeta^k$ $(k\in\zz)$.
Note that our terminology is slightly different from the existing literature, where
the simultaneous desingularization of our family is called the (Batyrev) mirror family.
We need to give a toric description of these hypersurfaces.
Let $\Delta^{\circ}$ be the polytope in the lattice $N\cong \zz^4\,$
with vertices, $v_1=(4,-1,-1,-1)$, $v_2=(-1,4,-1,-1)$,
$v_3=(-1,-1,4,-1)$, $v_4=(-1,-1,-1,4)$, $v_5=(-1,-1,-1,-1)$.
The fan $\Xi$ of the toric variety $\cpk/(\zz_5)^3$ is obtained by coning
over the faces of $\Delta^{\circ}$.
Then the $x_i$ can be regarded as the generators of the homogeneous coordinate
ring of $\cpk/(\zz_5)^3$ corresponding to the $v_i\,$.

\subsection{Orbifold Structure}

\begin{definition}
An orbifold structure on a Hausdorff, separable
topological space $X$ is given by an open cover $\mathbb{U}$ of $X$ satisfying
the following conditions:
\begin{itemize}
\item[(i)] Each element $U$ in $\mathbb{U}$ is uniformized, say by
$(V,G,\pi)$.  Namely, $V$ is a smooth manifold and $G$ is a finite
group acting smoothly on $V$ such that $U=V/G$ with $\pi$ as the
quotient map. Let ${\rm Ker}(G)$ be the subgroup of $G$ acting
trivially on $V$. \item[(ii)] For $U'\subset U$, there is a
collection of injections $(V',G',\pi') \to (V,G,\pi)$. Namely, the
inclusion $i: U' \subset U$ can be lifted to maps $\tilde{i}:V'
\to V$ and an injective homomorphism $i_{\ast}:G' \to G$ such that
$i_{\ast}$ is an isomorphism from ${\rm Ker}(G')$ to ${\rm
Ker}(G)$ and $\tilde{i}$ is $i_{\ast}$-equivariant. \item[(iii)]
For any point $x \in U_1 \cap U_2$, $U_1,U_2 \in \mathbb{U}$,
there is a $U_3 \in  \mathbb{U}$ such that $x \in U_3 \subset U_1
\cap U_2$.
\end{itemize}
\end{definition}

For any point $x \in X$, suppose that $(V,G,\pi)$ is a
uniformizing neighborhood and $\bar{x} \in \pi^{-1}(x)$. Let $G_x$
be the stabilizer of $G$ at $\bar{x}$. Up to conjugation, it is
independent of the choice of $\bar{x}$ and is called the
\emph{local group} of $x$.  Then there exists a sufficiently small neighborhood
$V_x$ of $\bar{x}$ such that $(V_x,G_x,\pi_x)$ uniformizes a small neighborhood of $x$,
where $\pi_x$ is the restriction\/ $\pi|_{V_x}$\/.  $(V_x,G_x,\pi_x)$ is called
a {\it local chart}\/ at $x$. The orbifold structure is called {\it reduced} if
the action of $G_x$ is effective for every $x$.

Let\/ $pr:E \to X$\/ be a rank $k$ complex {\it orbifold bundle}\/ over an orbifold $X$ 
(cf.$\;$\cite{[CR1]}).  
Then a uniformizing system for\/ $E|_U = pr^{-1}(E)$ over a uniformized subset 
$U$\/ of $X$\/ consists of the following data: 
\begin{itemize}
\item[(i)] A uniformizing system $(V,G,\pi)$ of $U$.
\item[(ii)] A uniformizing system $(V\times \cc^k,G,\tilde{\pi})$ for $E|_U$.
 The action of $G$ on $V\times \cc^k$ is an extension of the action of $G$ on 
 $V$ given by $g\cdot (x,v) =(g\cdot x, \rho (x,g)v)$ where 
$\rho:V\times G \to {\rm Aut}(\cc^k)$ is a smooth map satisfying: 
$$ \rho(g\cdot x,h) \circ \rho(x,g) = \rho(x,hg), \quad g,h \in G , \; x\in V.$$
\item[(iii)] The natural projection map $\tilde{pr}:V\times \cc^k \to V$ 
 satisfies $\pi \circ \tilde{pr} = pr \circ \tilde{\pi}$. 
\end{itemize}

By an orbifold connection $D$ on $E$ we mean an equivariant connection that satisfies $D=g^{-1}Dg$
for every uniformizing system of $E$. Such a connection can be always obtained by the 
averaging trick and an equivariant partition of unity.

\subsection{Twisted Sectors}

Let $X$\/ be an orbifold. Let $\widetilde{X}_k$  denote the
set of pairs $(x,({\bf g})_x)$ where $({\bf g})_x$ stands for the
conjugacy class of ${\bf g}=(g_1,\ldots\hspace{-1pt} ,g_k)$ with
$g_j\in G_x$\/.
Let $V_x^{g_j}$ denote the fixed point set of $g_j$ in $V_x$\/, and let
$C(g_j)$ denote the centralizer of $g_j$ in $G_x$\/.
$\widetilde{X}_k$ has a natural, possibly nonreduced, orbifold structure
 (cf.$\;$\cite{[CR1]},$\,$\cite{[Ka]}) whose local chart at
$(x,({\bf g})_x)$ is given by
\begin{equation}\label{X_k-orb-str}
\left(\, V_x^{{\bf g}}\: ,\: C({\bf g}) \: , \:\pi :  V_x^{{\bf g}}
\rightarrow V_x^{{\bf g}}/C({\bf g}) \, \right)
\end{equation}
where $V_x^{{\bf g}} = \bigcap_{j=1}^{k} V_x^{g_j}$, $C({\bf g})=
\bigcap_{j=1}^{k} C(g_j)$.
The topology on $\widetilde{X}_k$ is, of course, specified by declaring each
$V_x^{{\bf g}}/C({\bf g})$ to be an open set.

We describe the connected components of $\widetilde{X}_k$.
 Each point\/ $x\in X$ has a local chart
$(V_x,G_x,\pi_x)$ which gives a  uniformized
neighborhood $U_x = V_x/G_x$ of $x$.  If $y\in U_x$, up to conjugation there is a
unique monomorphism
$i_{\ast}: G_y \rightarrow G_x$ . For ${\bf
g}\in (G_y)^k$, the conjugacy class $i_{\ast}({\bf g})_y$ is well-defined. We
define an equivalence relation $i_{\ast}({\bf g})_y \cong ({\bf g})_y$.
Let $T_k$ denote the set of equivalence classes. By slight abuse
of notation, we use $({\bf g})$ to denote the equivalence class to
which $({\bf g})_x$ belongs.
We will usually denote an element of
$T_1$ by $(g)$. $\widetilde{X}_k$ is decomposed as a disjoint
union of connected components
$$ \widetilde{X}_k \: = \: \bigsqcup_{({\bf g})\in T_k}
X_{({\bf g})}\, ,
$$
where $X_{({\bf g})} =\left\{(x,({\bf g}')_x)\:|\:
 {\bf g}' \in (G_x)^{k} , \, ({\bf g}')_x \in
({\bf g}) \right\}$.
Note that for ${\bf g} = (1,\ldots\hspace{-1pt} ,1)$ we have
$X_{({\bf g})} \cong X$.
  A component $X_{({\bf g})}$ is called a {\it
twisted $k$-sector}, provided that ${\bf g}$ is not the identity.
$X_{(g)}$ is simply called a {\it twisted sector}. An almost complex,
complex or K\"ahler structure on $X$ induces an analogous structure on
$X_{({\bf g})}$ via (\ref{X_k-orb-str}).

Now define
\begin{align*}
T_3^0 &= \left\{
({\bf g})=(g_1,g_2,g_3)\in T_3 \:|\; g_1 g_2 g_3 = 1\:
\right\}, \\[5pt]
{\mathcal T}_3^0 &= \{ X_{({\bf g})} \:|\: ({\bf g})\in T_3^0 \}
\end{align*}

There is a one-to-one correspondence between $T_2$ and $T_3^0$ given by
$(g_1,g_2) \mapsto (g_1,g_2,(g_1g_2)^{-1})$.
We shall call an element of ${\mathcal T}_3^0$\/ a {\it tricyclic sector}\/.
The twisted sectors of Calabi-Yau hypersurfaces of a Fano toric
variety are described in \cite{[Po]}.
Following the same line of argument, we can identify
${\mathcal T}_3^0$ in the
present example.

\section{Tricyclic Sectors}
\subsection{${\mathcal T}_3^0$ for Simplicial Toric Varieties}\label{subsec:T_3^0-toric}

 Let $Y$\/ be the toric variety associated to a simplicial fan $\Xi$ in a
 $d$-dimensional lattice $N$.
 Denote the set of $n$-dimensional cones by $\Xi(n)$.
  For a cone $\tau \in \Xi$, denote the set of its primitive
 $1$-dimensional generators by $\tau[1]$, the corresponding affine open subset of
 $Y$\/
 by $\Ut$, and the corresponding torus orbit by $O_{\tau}$. We write $\nu \le \tau$ if the cone
 $\nu$ is a face of the cone $\tau$, and $\nu < \tau$ if it is a
 proper subface.
 $U_\tau = \bigsqcup_{\hspace{2pt}\nu \le \tau} O_\nu$.
 Let $M = {\rm Hom}(N,\zz)$ be the dual lattice of $N$ with dual pairing $\lan \: , \,\ran$.
 For any cone $\tau \in \Xi$,
 denote its dual cone in $M \otimes \rr$ by $\check{\tau}$.
 Let $S_\tau = \check{\tau} \cap M$. $\cc[S_\tau]$  is the $\cc$-algebra
 with generators $\chi^m$ for each
 $m\in S_\tau$ and relations $\chi^m\chi^{m'}=\chi^{m+m'}$.
 $U_\tau = {\rm Spec}( \cc[S_\tau])$. Define $R(\tau):= \{ \, \sum a_i f_i
 \; | \;
 f_i \in \tau[1],\: 0 \leq a_i < 1 \,\} \cap N$. We will describe the orbifold structure of
 $Y$.

 Let $\sigma$ be any $d$-dimensional cone of $\Xi$.
 The elements of $\sigma[1]$, $f_1,
 \ldots\hspace{-1pt} , f_d$, are linearly independent in $N \otimes \rr$.
 Let $\Ns$ be the sublattice of $N$ generated by $f_1, \ldots\hspace{-1pt} , f_d$.
 Let $\Gs:= N/\Ns$ be the quotient group. $\Gs$ is finite and Abelian.
 Let $\sigma'$ be the cone $\sigma$ regarded in $\Ns$. Let $\Ms$ be the
 dual lattice of $\Ns$ and let $\csigp$ be the
 dual cone of $\sigma'$ in $\Ms$.
 $\Usp= {\rm Spec}(\cc[\csigp \cap \Ms])$. Note that $\sigma'$ is a smooth cone in
 $\Ns$. So $\:\Usp \cong \cc^d$.

 There is a canonical dual pairing
 $\Ms /M \times N/\Ns \to \qq/\zz \to \cc^*$,
 the first map by the pairing $\lan \: , \,\ran$ and the second by
 $q \mapsto \exp (2\pi iq)$. Now $\Gs$ acts on $\cc[\Ms]$,
 the group ring of $\Ms$, by:
 $ n(\chi^m) = \exp(2\pi i\lan m,n\ran)\chi^m$, for $n \in N$ and $m \in \Ms$.
 Observe that $\Gs$ acts on $\Usp$ and $(\cc [\Ms])^{\Gs} = \cc [M]$.
 Thus $\Us=\Usp/\Gs$. Let $\pi_\sigma$ be the quotient map.
 So $\Us$ is uniformized by $(\Usp,\Gs, \pi_\sigma)$. For any $\tau<\sigma$,
 the orbifold structure on $\Ut$ is same as the one induced from the
 uniformizing system on $\Us$. Then $\{(\Usp,\Gs, \pi_\sigma) \: |\: \sigma \in \Xi(d)\}$
 defines a reduced orbifold structure on $Y$.

 It is convenient to have the following description of the local groups.
 Let $F$ be the nonsingular matrix with generators $f_1, \ldots\hspace{-1pt} ,f_d$ of
 $\sigma$ as rows.
 Then $\check{\sigma}'$ is generated in $\Ms$ by the column vectors
 $f^1,\ldots\hspace{-1pt} ,f^d$ of the matrix $F^{-1}$.
 So $\chi^{f^1}, \ldots\hspace{-1pt} ,\chi^{f^d}$ are the coordinates of $\Usp$.
 For any $r=(r_1, \ldots\hspace{-1pt} ,r_d)\in N$, the corresponding coset
 $[r]\in \Gs $
 acts on $\Usp$ in these coordinates as a
 diagonal matrix: diag$(e^{2\pi ic_1}, \ldots\hspace{-1pt} ,e^{2\pi ic_d})$ where
 $c_i = \lan r,f^i \ran$.
 Such a matrix is uniquely represented by a $d$-tuple $a=(a_1,\ldots\hspace{-1pt} ,a_d)$
 where $a_i \in [0,1)$ and $c_i=a_i+b_i$, $b_i \in \zz$. In matrix notation,
 $r F^{-1}=a + b \iff r = aF+bF.$
 We denote the integral
 vector $aF$ in $N$ by $r_a$ and the diagonal matrix corresponding
 to $a$ by $g_a$.  The correspondence
 ${g_a} \leftrightarrow {r_a}$ gives a bijection
  between the elements of $\Gs$ and the elements of $R(\sigma)$.

 Now we examine the orbifold chart
 induced by $(\Usp,\Gs,\pi_\sigma)$
 at any point $x\in \Us$. By the orbit decomposition, there is unique
 $\tau \le \sigma$ such that $x\in O_\tau$. Assume $\tau$ is generated
 by $f_1, \ldots\hspace{-1pt},f_j$, $ j\le d$. Let $z$ be a preimage of $x$
 with respect to $\pi_\sigma$.
 Then $\chi^{f^i}(z)=0\,$ iff $\,i \le j$.
 Let $\,\Gt:=\{  g_a \in \Gs \: | \;
 a_i=0  {\rm \  if \ }  j+1\le i\le d \} =\{\, g_a \in \Gs \: | \; r_a \in R(\tau) \}
 $.
 We can find a small neighborhood
 $W\subset (\cc^*)^{d-j}$ of $(z_{j+1}, \ldots\hspace{-1pt} ,z_d)
 $ such that the inclusions
 $\cc^j\times W\hookrightarrow \Usp$ and $\Gt \hookrightarrow \Gs$ induces an
 injection of uniformizing systems $(\cc^j\times W,\Gt,\pi) \hookrightarrow
 (\Usp,\Gs,\pi_\sigma)$ on some small open neighborhood $U_x$ of $x$. So
 we have $G_x = \Gt$ and an orbifold chart
 $(\cc^j \times W,\Gt,\pi)$.

  Now we determine ${\mathcal T}_3^0$.
  Take any $x \in Y$ with nontrivial local group. Then $x$ belongs to a
  unique $O_\tau$ such that $\tau$ is not the trivial cone.
  Pick any elements $g_a,g_b$ from $G_x = \Gt$. We shall find $Y_{(\bf{g})}$ where
  ${\bf{g}}= (g_a,g_b,(g_ag_b)^{-1})$.
  Let $\tau_a, \tau_b$
  be the faces of $\tau$, whose interiors contain $r_a$ and
  $r_b$ respectively.
  Let $\sigma$ be any $d$-dimensional cone containing $\tau$.
  Let $z$ be any
  point in $\Usp$. Suppose $z$ is fixed by both $g_a$ and $g_b$.
   Then $\chi^{f^i}(z)=0 $ whenever $f_i \in \tau_a \cup \tau_b$. Hence
  $\pi_{\sigma}(z) \in \overline{O}_{\tau_a} \cap \overline{O}_{\tau_b} 
  \cap \Us$. A local uniformizing system for $Y_{(\bf{g})}$ is given by
  $(V_x^{{\bf g}},G_x,\pi)$, where $$V_x^{{\bf g}}= \,(\cc^j \times W) \,\cap\,
  \{\chi^{f^i}=0, \;\: \forall f_i\in \tau_a \cup \tau_b \}.$$
  This leads us to observe that
  $\, \{(x,{\bf{g}}) \in Y_{(\bf{g})} |x \in \Us \}$  is complex analytically 
  isomorphic to $\overline{O}_{\tau_a} \cap \overline{O}_{\tau_b}\cap \Us$.
  Since this is true irrespective of the choice of $\sigma$, 
  $Y_{\bf{g}} \cong \overline{O}_{\tau_a} \cap \overline{O}_{\tau_b}$.
  Note that $\overline{O}_{\tau_1} \cap \overline{O}_{\tau_2}$ is empty whenever
  $\tau_1[1] \cup \tau_2[1]$ does not generate an element of $\Xi$.

\begin{proposition} If $\tau_1[1] \cup \tau_2[1]$ generate an element of $\Xi$,
  then for every pair $r_1 \in R(\tau_1)\cap {\rm Int}(\tau_1)$,
 $r_2 \in R(\tau_2)\cap {\rm Int}(\tau_2)$
  we have a unique element of ${\mathcal T}_3^0(Y)$, which is analytically 
 isomorphic to $\overline{O}_{\tau_1} \cap \overline{O}_{\tau_2}$.
 As we vary over $\tau_1, \tau_2$, we obtain all elements of ${\mathcal T}_3^0(Y)$.
\end{proposition}

\subsection{${\mathcal T}_3^0$ for a Nondegenerate Quasi-smooth Hypersurface}
\label{subsec:T_3^0Hyper}

 A hypersurface $X$ of the toric variety $Y$ is called {\it quasi-smooth} if for
 any $\sigma \in \Xi(d)$, $\pi_{\sigma}^{-1}(X \cap \Us)$ is smooth. Then $X$ is clearly
 a suborbifold of $Y$, the orbifold structure being induced by
 $\{(\pi_{\sigma}^{-1}(X \cap \Us),\Gs, \pi_\sigma) \: |\: \sigma \in \Xi(d)\}$.

 $X$ is called {\it nondegenerate} if for every $\tau \in \Xi$,
 $X \cap O_{\tau}$ is either empty or a smooth variety of codimension one in $O_{\tau}$.
 When X is nondegenerate, the above orbifold structure on $X$ is reduced. So any point
 $x \in X \cap O_{\tau}$ has local group $G_\tau$. Then by an argument analogous to
 the  one for ${\mathcal T}_3^0(Y)$, we have the following description of
${\mathcal T}_3^0(X)$.

\begin{proposition} Let $X$ be a nondegenerate quasi-smooth hypersurface  
 of a simplicial toric variety $Y$ with fan $\Xi$.
 If $\tau_1[1] \cup \tau_2[1]$ generate an element of $\Xi$,
  then for every pair $r_1 \in R(\tau_1)\cap {\rm Int}(\tau_1)$,
 $r_2 \in R(\tau_2)\cap {\rm Int}(\tau_2)$
  the connected components of
  $X \cap \overline{O}_{\tau_1} \cap \overline{O}_{\tau_2}$
  are elements of ${\mathcal T}_3^0(X)$.   As we vary over $\tau_1, \tau_2$,
  we obtain all elements of ${\mathcal T}_3^0(X)$.
\end{proposition}

\subsection{$\mathcal{T}_3^0$ for Mirror Quintic}

 Let $Y = \cpk/(\zz_5)^3 $ and let $X$\/ be a
 member of the mirror quintic family such that $\psi \neq -5\zeta^k $
(\ref{eq:quintic}). Then $X$\/ is nondegenerate and quasi-smooth.
 The fan $\Xi$ of $Y$ is described in Subsection~\ref{subsec:quintic}.
 We shall conveniently use the action of an element of $\Gt$ on the 
 homogeneous coordinate  ring of $Y$ to describe that element. 
 For instance, $(\zeta,1,\zeta^2,\zeta^2,1)$ will
 represent an element of $\Gt$ where $\tau[1]=\{v_1,v_3,v_4\}$.
 The corresponding element
 of $R_\tau$ is $(\frac{1}{5}v_1 + \frac{2}{5}v_3 + \frac{2}{5}v_4)$.

 Note that $\bar{O}_\tau = \{x_i=0 \: | \: i \, 
 {\rm \ such \  that\ } v_i \in \tau[1] \}$.
 Since there is no complex reflection, if $x \in O_\tau $ has 
 nontrivial local group then
 $\dim(\tau) \ge 2$ and consequently ${\rm codim}(O_\tau) \ge 2 $.
 By nondegeneracy, $\dim(X \cap \bar{O}_{\tau}) = {\rm codim}(\tau) -1 $.
 Hence a nontrivial element of
 ${\mathcal T}_3^0(X)$ (i.e. an element not corresponding to 
 $({\bf g}) = (1,1,1) \, $) has dimension less than $2$. We have the 
 following classification of nontrivial elements of $\mathcal{T}_3^0(X)$.

\begin{corollary}\label{cor:cases}
 $X_{({\bf g})}$ is\/ $0$-dimensional
 if\/ $\tau_1[1] \cup \tau_2[1]$ generate a\/
 $3$-dimensional cone, i.e. if any of the following holds:
\begin{itemize}
\item[(i)]
 $\tau_1$,$\tau_2$ are both $2$-dimensional and have a $1$-dimensional face in common.
 \item[(ii)] $\tau_1$ is a $3$-dimensional and $\tau_2$ is a $2$-dimensional
 subface of $\tau_1$, or vice versa.
 \item[(iii)] $\tau_1=\tau_2$ and\/ $\dim(\tau_i)=3$.
\end{itemize}
 $X_{({\bf g})}$ is\/ $1$-dimensional if\/ $\tau_1[1] \cup \tau_2[1]$ generate a\/
 $2$-dimensional cone, i.e. if\/ $\tau_1 = \tau_2$ and\/ $\dim(\tau_i) = 2$.
 \end{corollary}

 From the description of $\Gt$ in terms of $R_{\tau}$
 in Subsection~\ref{subsec:T_3^0-toric}, it is not hard to see that
 $\Gt \cong (\zz_5)^{{\rm dim}(\tau)-1}$. For a point $(x,({\bf g})) \in X_{({\bf g})}$
 with $x \in O_{\tau} \cap X$ and $({\bf g})$ determined by 
 $r_i \in R(\tau_i)\cap{\rm Int}
 (\tau_i)$, a local uniformizing system for $X_{({\bf g})}$ is given by 
 $(V_x^{\bf{g}}, \Gt , \pi )$, where $V_x^{\bf{g}}$ is a small neighborhood
 of a preimage of $x$ in $\widehat{X} \cap \{x_i=0, \; \forall i 
 {\rm\   such \  that\ } v_i \in \tau_1 \cup \tau_2 \}$.
 Here $\widehat{X}\subset\cpk$ is the zero locus of the polynomial (\ref{eq:quintic}).

\section{Chen-Ruan Orbifold Cohomology}\label{sec:CR cohomology}

Assume that $X$\/ is a $d$-dimensional compact almost complex orbifold with
 an almost complex
 structure $J$\/ (cf.$\;$\cite{[CR1]}).
Then for a point $x$ with nontrivial local group
$G_x$,
$J$\/
 gives rise to an effective representation $\rho_x:G_x\to GL(d,\cc)$. For any
 $g\in G_x$ we write $\rho_x(g)$, up to conjugation, as a diagonal matrix
 $${\rm diag}(e^{2\pi i\frac{m_{1,g}}{m_g}}, \,
\ldots ,e^{2\pi i\frac{m_{d,g}}{m_g}}),$$
 where $m_g$ is the order of $g$ in $G_x$, and $0\leq m_{i,g}<m_g \, .$  Define
 a function $\iota:\widetilde{X}_1 \to \qq$ by
$$ \iota(x,(g)_{x}) \:=\left\{ \begin{array}{ll}
\sum_{i=1}^d \frac{m_{i,g}}{m_g} \hspace{5pt} &{\rm if \ }g\neq 1 \, , \\[5pt]
0   &  {\rm if \  }g = 1\, .
\end{array}\right.$$
 This function $\iota:\widetilde{X}_1 \to \qq$\/ is locally constant.
 Denote its
 value on $X_{(g)}$ by $\iota_{(g)}$.  We call $\iota_{(g)}$ the
 \emph{degree shifting number} of $X_{(g)}$\/.
 It has the following properties:
\begin{itemize}
\item[(i)] $\iota_{(g)}$ is integral iff\/ $\rho_x(g)\in
SL(d,\cc)$. \item[(ii)] $\iota_{(g)} + \iota_{(g^{-1})} = {\rm
rank} (\rho_x(g)-Id)=d-\dim_{\cc}X_{(g)}$.
\end{itemize}

A $C^\infty$ differential form on $X$ is a $G$-invariant $C^\infty$ 
differential form on $V$ for each uniformizing system $(V,G,\pi)$. 
Then {\it orbifold integration} is defined as follows. Suppose $U=V/G$
is connected. For any compactly supported differential $d$-form 
$\omega$ on $U$, which is, by definition, a $G$-invariant $d$-form
$\tilde{\omega}$ on $V$, 
\begin{equation}\label{eq:orb-integral}
 \int_U^{orb} \omega := \frac{1}{|G|} \int_V
 \tilde{\omega},
\end{equation}
 where $|G|$ is the order of $G$. Then orbifold
integration over $X$ is defined by using a $C^\infty$ partition of
unity. The orbifold integration coincides with the usual measure 
theoretic integration if and only the orbifold structure on $X$   
is reduced. 

Holomorphic forms for a complex orbifold $X$ are again obtained by
patching $G$-invariant holomorhic forms on the uniformizing systems
$(V,G,\pi)$. We consider the \v{C}ech cohomology groups of X
and $X_{\bf{g}}$ with coeffcients in the sheaves of holomorphic forms.
These \v{C}ech cohomology groups can be identified with the Dolbeault
cohomology groups of $(p,q)$-forms \cite{[Bai]}.

 \begin{definition}(cf.$\,$\cite{[CR1]}) Let $X$ be a closed complex orbifold. 
We define the \emph{orbifold cohomology groups}\/ of $X$\/ by 
$$ H^n_{orb}(X) := \bigoplus_{(g)\in T_1 } H^{n-2\iota_{(g)}}(X_{(g)};\cc). $$
We
 define, for $0 \le p,q \le \dim_{\cc}X,$ \emph{orbifold Dolbeault cohomology groups}
 $$H_{orb}^{p,q}(X) := \bigoplus_{(g)\in T_1 }H^{p-\iota_{(g)},q-\iota_{(g)}}(X_{(g)};
\cc).$$
\end{definition}

\subsection{Obstruction Bundle}

Choose $({\bf g})=(g_1,g_2,g_3)\in T_3^0$. Let $(x,({\bf g})_x)$
be a generic point in $X_{({\bf g})}$. Let $\kG$\/ be the subgroup
of $G_x$ generated by $g_1$ and $g_2$. Consider an orbifold
Riemann sphere with three orbifold points, $(\sss^2,
(x_1,x_2,x_3),$\linebreak[1]$(k_1,k_2,k_3))$.  When there is no
confusion, we will simply denote it by $S^2$.  The orbifold
fundamental group is
$$ \pi_1^{orb} (S^2) = \left\{ \: \lambda_1 ,
\lambda_2 , \lambda_3 \;|\;  \lambda_i^{k_i}=1, \: \lambda_1
\lambda_2 \lambda_3 =1  \: \right\},  $$
where $\lambda_i$ is represented by a loop around the marked $x_i$.
There is a surjective homomorphism
\begin{equation}\label{map:rho}
 \rho : \pi_1^{orb} (S^2) \rightarrow \kG \, ,
\end{equation}
specified by mapping $\lambda_i \mapsto g_i$. Ker($\rho$) is a
finite-index subgroup of $\pi_1^{orb} (S^2)$. Let $\tilde{\Sigma}$
be the orbifold universal cover of $S^2$. Let $\Sigma =
\tilde{\Sigma}/{\rm Ker}(\rho)$.  Then $\Sigma$ is smooth, compact
and $\Sigma/\kG=S^2$.  The genus of $\Sigma$ can be computed using
Riemann-Hurwitz formula for Euler characteristic of a branched
cover, and turns out to be
\begin{equation}\label{eq:genus}
g(\Sigma) \:=\: \frac{1}{2}\left(2+|\kG|-\sum_{i=1}^{3} \frac{|\kG|}{k_i}
\right) \, .
\end{equation}
$\kG$ acts holomorphically on $\Sigma$ and hence $\kG$ acts on $H^{0,1}(\Sigma)$.
The ``obstruction bundle'' $E_{({\bf g})}$ over $X_{({\bf g})}$ is constructed 
as follows. On the local chart $(V_x^{{\bf g}}\: ,\: C({\bf g}) \: , \:\pi)$
of $X_{({\bf g})}$, $E_{({\bf g})}$  is given by 
$(TV_x \otimes H^{0,1}(\Sigma))^{\kG} \times V_x^{{\bf g}} \to V_x^{{\bf g}}$,
where $(TV_x \otimes H^{0,1}(\Sigma))^{\kG}$ is the $\kG$-invariant subspace. We 
define an action of $C(\bf{g})$ on $TV_x \otimes H^{0,1}(\Sigma)$, which is
the usual one on $TV_x$ and trivial on $H^{0,1}(\Sigma)$. Then the actions of
$C(\bf{g})$ and $\kG$ commute and $(TV_x \otimes H^{0,1}(\Sigma))^{\kG} $ is invariant
under $C(\bf{g})$. Thus we have obtained an action of $C(\bf{g})$ on 
$(TV_x \otimes H^{0,1}(\Sigma))^{\kG} \times V_x^{{\bf g}} \to V_x^{{\bf g}}$, 
extending the usual one on $V_x^{{\bf g}}$. These trivializations fit together
to define the bundle $E_{({\bf g})}$ over $X_{({\bf g})}$. If we set 
 $e: X_{({\bf g})} \rightarrow X$ to be the  map
given by $(x,({\bf g})_x)\mapsto x$\/, one may think of $E_{({\bf g})}$ 
as $ ( e^{\ast}TX \otimes H^{0,1}(\Sigma))^{\kG} \,$.
The rank of $E_{({\bf g})}$ is given by the formula
\begin{equation}\label{eq:dimension-Eg}
{\rm rank}_{\cc}(E_{({\bf g})}) \:=\: \dim_{\cc}(X_{({\bf g})}) -
 \dim_{\cc} (X) \,+\, \sum_{j=1}^{3} \iota_{(g_j)}
\, .
\end{equation}

\subsection{Orbifold Cup Product}
There is a natural map $\Phi : X_{(g)} \rightarrow X_{(g^{-1})}$ defined by 
$(x,(g)_x)\mapsto (x,(g^{-1})_x)$.  
\begin{definition}\label{def:2-point}
Let $d=\dim_{\cc}(X)$.  For any integer $0\leq n \leq 2d$, the pairing 
$$
\langle \; ,\: \rangle_{orb} \: :\: H^n_{orb}(X) \times H^{2d-n}_{orb}(X)\longrightarrow \cc
$$
is defined by taking the direct sum of 
$$
\langle \; ,\: \rangle_{orb}^{(g)} \: : \: H^{n-2\iota_{(g)}}(X_{(g)};\cc) \times 
H^{2d-n-2\iota_{(g^{-1})} }(X_{(g^{-1})};\cc) \longrightarrow \cc
$$
where
$$
\langle \alpha , \beta\rangle_{orb}^{(g)} \:=\: \int^{orb}_{X_{(g)}} 
\alpha\wedge\Phi^{\ast}(\beta)
$$
for $\alpha \in H^{n-2\iota_{(g)}}(X_{(g)};\cc)$, and
$\beta\in H^{2d-n-2\iota_{(g^{-1})} }(X_{(g^{-1})};\cc)$. 
\end{definition}

Choose an orbifold connection $A$\/ on $E_{({\bf g})}$.
Let $e_A(E_{({\bf g})})$ be the Euler form computed from the connection
$A$\/ by Chern-Weil theory.
Let $\eta_j\in H^{d_j}(X_{(g_j)};\cc)$,  for\/ $j=1,2,3$.
Define maps $e_j: X_{({\bf g})}\rightarrow X_{(g_j)}$\/ by
$(x,({\bf g})_x) \mapsto (x,(g_j)_x)$.
\begin{definition}
We define the $3$-{\it point function} to be
\begin{equation}\label{eq:3point}
 \langle \eta_1 , \eta_2 , \eta_3 \rangle_{orb} \:= \:
\int_{X_{({\bf g})}}^{orb} e_1^{\ast} \eta_1 \wedge e_2^{\ast}\eta_2
\wedge e_3^{\ast}\eta_3 \wedge
e_A(E_{({\bf g})}) \, .
\end{equation}
Note that the above integral does not depend on the choice of $A$.
As in Definition~\ref{def:2-point}, we extend the $3$-point function to 
$H^{\ast}_{orb}(X)$\/ via linearity.  
We define the orbifold cup product by the relation
\begin{equation}\label{eq:cup-product}
\langle\eta_1 \cup_{orb} \eta_2\, ,\eta_3 \rangle_{orb}\:=\:
\langle \eta_1 , \eta_2 , \eta_3 \rangle_{orb} \, .
\end{equation}
Again we extend $\,\cup_{orb}$ to $H^{\ast}_{orb}(X)$\/ via linearity.
\end{definition}

\begin{remark}
Note that if $({\bf g})=(1,1,1)$, then\/ $\eta_1 \cup_{orb}\eta_2$
is just the ordinary cup product\/ $\eta_1 \cup \eta_2$ in
$H^{\ast}(X)$.
\end{remark}

\smallskip

\section{Computations for Mirror Quintic}\label{sec:computations}

 Let $X$ be a generic member of the mirror quintic family.
The crucial part of the cup product computation is identifying
the obstruction bundles $E_{({\bf g})}$ and integrating their Euler
class. We will content ourselves with computing the $3$-point 
functions.  We
showed in Subsection~\ref{subsec:T_3^0Hyper} that $X_{({\bf g})}$
is either a point or a curve, provided that $({\bf g})$ is not the
identity.  The ordinary cup product corresponding to the identity
case will be covered in Section~\ref{sec:ord-cup}.

\subsection{Point Case}
If $X_{({\bf g})}$ is a point then $E_{({\bf g})}$ is a vector
space and we can use rank formula (\ref{eq:dimension-Eg}).
For the 3-point function to be nonzero, we need ${\rm rank}_{\cc}
(E_{({\bf g})})=0$, which forces  
$\sum_{j=1}^{3} \iota_{(g_j)} = 3$.  In this case, we have  
$\langle \eta_1, \eta_2, \eta_3 \rangle_{orb} = \frac{1}{25}\,
\eta_1\eta_1\eta_3$, where
$\eta_j \in H^{0,0}(X_{(g_j)})$.  

All in all there are 930 possible choices of point sectors ({\bf g}) that give 
nonzero 3-point function.   We can divide these into two types: 
\begin{itemize}
\item[(i)]  $\iota_{(g_j)}=1\,$ for all\/ $j=1,2,3$.  
\item[(ii)]  $\iota_{(g_3)}=0\,$, i.e.\/ $g_3 = {\rm id}$. 
\end{itemize}
Note that for type (i), $\iota_{(g_1 g_2)} = 2\,$ by property (ii) of degree shifting 
number in Section~\ref{sec:CR cohomology}.  There are 
810 point sectors of type (i) and 120 point sectors of type (ii).  
For example, 
\begin{align*}
((\zeta^4,\zeta,1,1,1)\, ,\,(1,\zeta^3,\zeta^2,1,1)\, ,\,(\zeta,\zeta,\zeta^3,1,1))\, ,\\
((\zeta^3,\zeta,\zeta,1,1)\, ,\,(\zeta,\zeta^2,\zeta^2,1,1)\, ,\,
(\zeta,\zeta^2,\zeta^2,1,1))\, ,\\
((\zeta^2,\zeta^3,1,1,1)\, ,\,(\zeta^2,\zeta,\zeta^2,1,1)\, ,\,(\zeta,\zeta,\zeta^3,1,1))\,\;
\end{align*}
are type (i) and 
$$((\zeta^3,\zeta,\zeta,1,1)\, ,\,(\zeta^2,\zeta^4,\zeta^4,1,1)\, ,\,(1,1,1,1,1))$$
is type (ii).  

\smallskip

\subsection{Curve Case}

Now consider the case when $X_{({\bf g})}$ is a curve. Then by
Corollary~\ref{cor:cases}, up to a canonical isomorphism ,
$X_{({\bf g})} = \bar{O}_{\tau} \cap X$, where $\dim(\tau)=2$.
Consequently $\Gt \cong \zz_5$. Hence $\kG$, the subgroup
of $\Gt$ generated by $g_j$ ($j=1,2,3$), is just $\Gt$.
Thus we may have $g_2 = (g_1)^k$ where $k=1,2,3,4$. The
kernel of the homomorphism $\rho$ (\ref{map:rho}), and
hence the Riemann surface $\Sigma$, depend
on the value of $k$. We shall fully compute one representative
case for each value of $k$. The computations for remaining cases can be
completed using the same method and are left to the reader.

In the case $k=4$, $g_3 = 1$, and so to construct $E_{({\bf g})}$
we have to consider the orbifold sphere
$(\sss^2,(x_1,x_2,x_3),(5,5,1))$. Then $\Sigma = \sss^2$, the
smooth sphere.

\begin{figure}[!ht]
\begin{center}
\includegraphics[scale=.9]{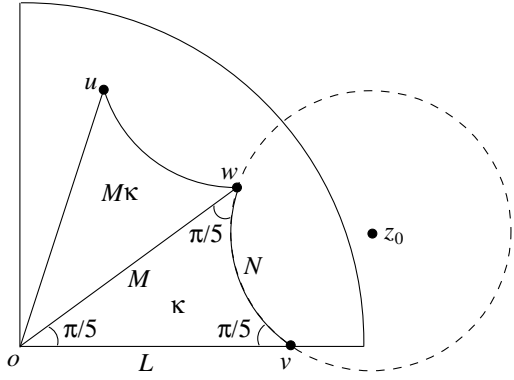}
\end{center}
\caption{Fundamental region $\epsilon$ in the first quadrant}
\label{fig:s2}
\end{figure}

 For other values of $k$, the
 order of each $g_j$ is $5$. So in these cases we need to work with
 the orbifold sphere $(\sss^2,(x_1,x_2,x_3),(5,5,5))$,
 which we continue to denote by $S^2$.  Recall that
$$ \pi_1^{orb} (S^2) = \left\{ \: \lambda_1 ,
\lambda_2 , \lambda_3 \;|\;  \lambda_i^{5}=1, \: \lambda_1
\lambda_2 \lambda_3 =1  \: \right\}. $$ Its orbifold universal
cover is the hyperbolic plane $\hh^2$ (cf.$\,$\cite{[Sc]}). We
will use the Poincar\'e disk model of $\hh^2$. Let $\kappa$ be the
geodesic triangle $\triangle ovw$ in $\hh^2$ as in
Figure~\ref{fig:s2} with all of its angles equal to $\theta
=\pi/5$. The sides of $\kappa$ are labelled by $L,M,N$. Without
causing too much confusion, we will denote the hyperbolic
reflection about a side of $\kappa$ by the same letter.  Let
$\epsilon$ denote the region $\kappa \cup M\kappa$.

There is a properly discontinuous action of $\pi_1^{orb} (S^2)$ on
$\hh^2$ which gives a tiling of $\hh^2$ with fundamental region
$\epsilon$\/ as follows: $\lambda_1$ acts as the composition
$L\circ M$, $\lambda_2$ acts as the composition $M\circ N$, and
$\lambda_3$ acts as the composition $N\circ L$. The quotient of
$\hh^2$ by $\pi_1^{orb} (S^2)$ action is the fundamental region
$\epsilon$ with the edge $\overline{ou}$ identified with the edge
$\overline{ov}$, and the edge $\overline{uw}$ identified with the
edge $\overline{vw}$.

It is not difficult to compute that
$$v \:=\:
\sqrt{\frac{2}{1+\sqrt{5}}} \:
\approx \: 0.786151 \; .$$
The side $N$\/ of $\kappa$ is an arc of the circle in $\cc$\/ with
center $z_0$ and radius $R$, where
\begin{align*}
z_0 \: &=\:  \sqrt{\frac{5+3\sqrt{5}}{10}} 
\exp(\pi i /10)
\:\approx\: 1.02909 + 0.33437\,i
\, , \\[5pt]
R \:& =\:  {\rm Im}\,z_0\cdot\sec\theta  \: =\:  \frac{-1+\sqrt{5}}{2\sqrt[4]{5}} \:
\approx \: 0.413304 \: .
\end{align*}
Then we can write $\lambda_1(z)=e^{-i4\theta}z\,$,\/ and
$$\lambda_2(z) = e^{i4\theta}\cdot
\frac{\overline{z_0}\, z + (R^2 - |z_0|^2)}{z-z_0}\, .$$

\subsubsection{The Case $g_1 = g_2$}

\begin{figure}[!ht]
\begin{center}
\includegraphics[scale=.65]{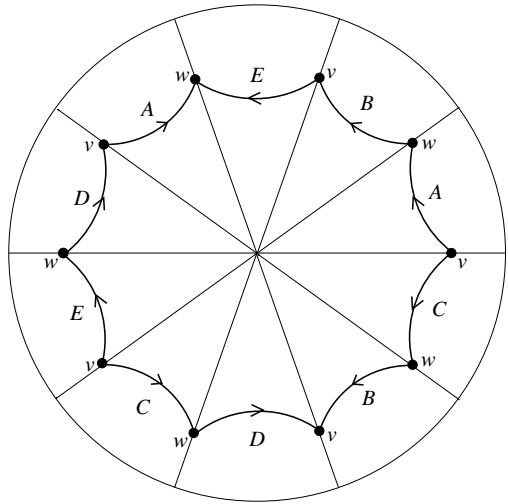}
\end{center}
\caption{Genus two surface $\Sigma$}
\label{fig:tengon}
\end{figure}

\hspace{2pt}Consider the homomorphism $\rho: \pi_1^{orb} (S^2)
\rightarrow \kG\,$ given by $\lambda_i \mapsto g_i$. ${\rm
Ker}(\rho)$ is a normal subgroup generated by commutators of
$\lambda_i$ and the element $\lambda_1 \lambda_2^{-1}$. ${\rm
Ker}(\rho)$ acts freely on $\hh^2$ with quotient being the decagon
shown in Figure~\ref{fig:tengon}. The sides of the decagon are
identified by the elements of ${\rm Ker}(\rho)$ given in
Table~\ref{table:tengon}.
\begin{table}[!ht]
\caption{ }
\begin{center}\label{table:tengon}
\begin{tabular}{|c|c|c|c|c|}
\hline
\vspace{-.3cm}
&&&&\\
     $A$  & $B$ & $C$ & $D$ & $E$\\ \vspace{-.3cm}
&&&&\\ \hline
\vspace{-.3cm}
&&&&\\ \vspace{-.3cm}
\hspace{10pt}$\lambda_1^{-1}\lambda_2$\hspace{10pt}  &
\hspace{10pt}$\lambda_1\lambda_2^{-1}$\hspace{10pt} &
$\lambda_1(\lambda_1\lambda_2^{-1})\lambda_1^{-1}$
& $\lambda_1^2(\lambda_1\lambda_2^{-1})\lambda_1^{-2} $
& $\lambda_1^3(\lambda_1\lambda_2^{-1})\lambda_1^{-3} $ \\
&&&&\\ \hline
\end{tabular}
\end{center}
\end{table}

\begin{figure}[!ht]
\begin{center}
\includegraphics[scale=.9]{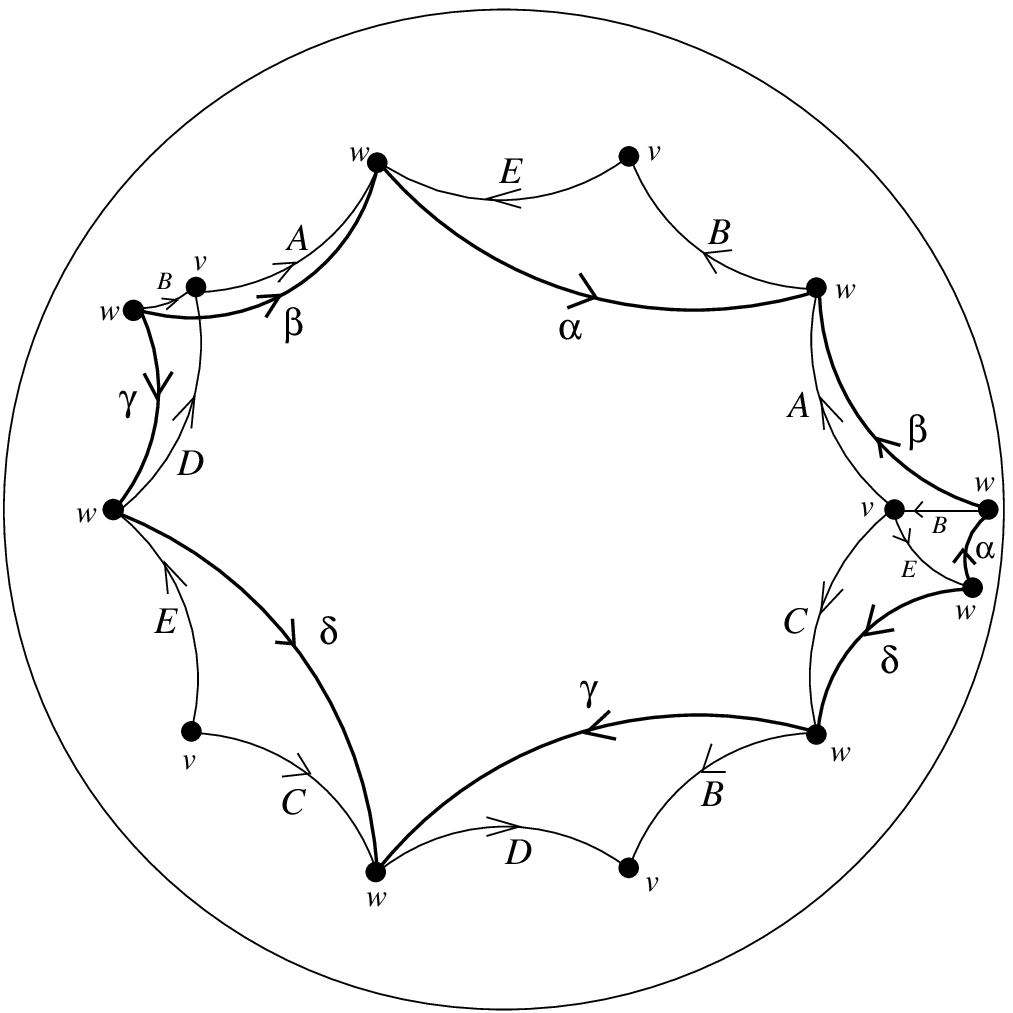}
\end{center}
\caption{ }
\label{fig:octagon}
\end{figure}

We cut out four hyperbolic triangles along the geodesic curves,
$\alpha = (BE)^{-1}, \,
\beta = BA , \, \gamma = BD^{-1} , \, \delta = E^{-1} C$.  After pasting the
triangles
along the geodesic curves, $A,B,C,D,E$\/ or their translates by elements of
${\rm Ker}(\rho)$, we obtain
the octagon in Figure~\ref{fig:octagon}, whose boundary
curves comprise a symplectic basis
$\left\{a_1, b_1 , a_2, b_2 \right\}=
\left\{\alpha, \beta, \gamma,
\delta \right\}$  
for $H_1(\Sigma ;\zz)$\/ after appropriate
identifications.

Taking the Universal Coefficients Theorem duals,
we obtain a canonical symplectic basis  
$\left\{\check{\alpha},  \check{\beta}, \check{\gamma},
\check{\delta} \right\}$ for $H^1(\Sigma;\zz)$.
From the Riemann bilinear relations in
\cite{[Gr]} (Chapter$\;$V, \S 3), we know that there exists a complex basis
$\left\{ \omega_1, \omega_2 \right\}$ for $H^{1,0}(\Sigma)
= H^0(\Sigma , K_{\Sigma})$\/  such that
its period matrix with respect to the $a_i$-classes is the identity matrix,
and its period matrix with respect to the $b_i$-classes is a symmetric
complex matrix
$$
P=\left(\begin{array}{cc}p&q \\q&s \end{array}\right)
$$
with the imaginary part, ${\rm Im} P$, being a positive definite matrix.
Hence we can write $\{ \omega_1 , \omega_2\} = \{
\check{\alpha}+p\check{\beta}+q\check{\delta},
\check{\gamma}+q\check{\beta}+s\check{\delta}
\}$.

The rotation $\lambda_1$ induces the following map on $\pi_1(\Sigma)$:
$$
\begin{array}{l}
\alpha \:\longmapsto\: A^{-1}C = \beta^{-1} \alpha^{-1} \delta \\[7pt]
\beta \:\longmapsto\: C^{-1}B^{-1} = \delta^{-1} \alpha \\[7pt]
\gamma \:\longmapsto\: C^{-1}E = \delta^{-1} \\[7pt]
\delta \:\longmapsto\: A^{-1}D^{-1} = \beta^{-1} \gamma
\end{array}
$$
Hence
the automorphism $(\lambda_1)_{\ast}:H_1(\Sigma)\rightarrow
H_1(\Sigma)$ can be expressed in the matrix
\begin{equation}\label{matrix:H_1}
\left(\begin{array}{rrrr}-1&1&0&0 \\
-1&0&0&-1 \\
0&0&0&1 \\
1&-1&-1&0
\end{array}\right)
\end{equation}
with respect to the basis\/ $\{ \alpha, \beta, \gamma, \delta \}$.
The automorphism $\lambda_1^{\ast}:H^1(\Sigma)\rightarrow
H^1(\Sigma)$ can be expressed in the matrix
\begin{equation}\label{matrix:H^1}
\left(\begin{array}{rrrr}-1&-1&0&1 \\
1&0&0&-1 \\
0&0&0&-1 \\
0&-1&1&0
\end{array}\right)
\end{equation}
with respect to the dual basis $\left\{\check{\alpha},  \check{\beta}, \check{\gamma},
\check{\delta} \right\}$, which is the transpose of matrix (\ref{matrix:H_1}).
From matrix (\ref{matrix:H^1}), we can easily calculate that
$$
\begin{array}{l}
\lambda_1^{\ast} \omega_1 = (-1-p+q) \check{\alpha} +(1-q)\check{\beta} +(-q)\check{\gamma}
+(-p)\check{\delta}\, ,
\\[7pt]
\lambda_1^{\ast}\omega_2 =  (-q+s) \check{\alpha} +(-s)\check{\beta} +(-s)\check{\gamma}
+(1-q)\check{\delta}\, .
\end{array}
$$

Since $\lambda_1$ is a holomorphic map, $\lambda_1^{\ast}$
preserves the subspace $H^{1,0}(\Sigma)$.  Hence
$\lambda_1^{\ast}\omega_1$ and $\lambda_1^{\ast}\omega_2$ can be
expressed as complex linear combinations of $\omega_1$ and
$\omega_2$, so we easily obtain a system of equations satisfied by
the triple $(p,q,s)$:
\begin{equation}\label{eq:pqs}
\left\{ \begin{array}{l}
1-q = p(-1-p+q)-q^2 \\[7pt]
-p= q(-1-p+q) - qs \\[7pt]
-s= p(-q+s) -sq\\[7pt]
1-q=q(-q+s) -s^2
\end{array}\right.
\end{equation}
One can verify that
there are four solutions to system~(\ref{eq:pqs}), but only one of them
satisfies the condition that Im$P$ is positive definite, namely,
\begin{equation}\label{pqs}
\left\{\begin{array}{l}
p=\exp(3\pi i /5) \; ,
\\[7pt]
q=\frac{1}{2}\left( 1+i\sqrt{5-2\sqrt{5}}\right), \\[7pt]
s=\exp(2\pi i /5)=\zeta \; .
\end{array}
\right.
\end{equation}
We now know that $H^{0,1}(\Sigma)=\overline{H^{1,0}(\Sigma)}$\/ is generated by
$$\{\overline{\omega}_1,
\overline{\omega}_2 \} = \left\{ \check{\alpha} + \overline{p}\check{\beta} +
\overline{q}\check{\delta}, \check{\gamma} + \overline{q}\check{\beta}
+\overline{s}\check{\delta} \right\}\, .$$

\subsubsection*{Representative Case}

 Choose\/ $g_1 = (\zeta^n , \zeta^{5-n} , 1, 1, 1)$\/ where\/ $n=1,2,3$ or 4.  
 First we consider $e^{\ast}TX$
 over the corresponding $X_{({\bf g})}$, which we identify with $X\cap \{x_1 =x_2=0 \}$
 as an analytic space. The orbifold structure on $X_{({\bf g})}$ is not reduced. At the
 three points\/ $p_j=X\cap \{x_1 =x_2=x_j=0 \}$\/ corresponding to\/ $j=3,4$ or 5, 
 the local groups are isomorphic to $\zz_5^2$ and contain $\kG \cong \zz_5$ as a subgroup.
 The local group is $\kG$
 at the remaining points. But the $\kG$ action in the local chart is trivial at each 
 point. By quotienting out all the local groups by $\kG$, we can associate a reduced
 orbifold $X_{(\bf{g})}':=X_{({\bf g})}^{red}$ 
to $X_{({\bf g})}$. We again identify $X_{(\bf{g})}'$ 
 with $X\cap \{x_1 =x_2=0 \}$, which is a genus zero curve 
(cf.$\;$\cite{[Po]},$\;$Section 5.3).

Take the open charts $V_i=\{x_i\neq 0\}$ of $\cpk$, and $U_i=V_i/(\zz_5)^3$ of
$\cpk/(\zz_5)^3$,
 where $x_i$ are the homogeneous coordinates.  Let $\pi_i :V_i \rightarrow U_i$ denote
the projection map.
 For points in $V_5$, we define new coordinates $z_i = {x_i}/{x_5}$, with $z_5=1$.
On $V_4$ we define coordinates $w_i=x_i/x_4$, with $w_4=1$.
Note that
$$  z_i = \frac{x_i}{x_5} = \frac{x_i/x_4}{x_5/x_4} =
\frac{w_i}{w_5} \, ;  \quad
w_i  = \frac{x_i}{x_4} = \frac{x_i/x_5}{x_4/x_5}= \frac{z_i}{z_4}\, .
$$
Let $\widehat{X}_{i} = \pi_i^{-1}(X \cap U_i)$ and let 
${X}_{({\bf g})}^{\circ} =X_{({\bf g})}\setminus \{ p_3,p_4,p_5\}$, 
$X_{({\bf g})}^{'\circ} =X_{({\bf g})}'\setminus \{ p_3,p_4,p_5\}$.

Next we write down a local framing for $
T\widehat{X}_{5}\bigr|_{\pi_5^{-1}(X_{({\bf g})}^{\circ})}$ in terms of
the $z_i$ coordinates.  Note that
$$
\widehat{X}_5 \:=\: \{ z_1^5 +z_2^5 +z_3^5+ z_4^5 +1 + \psi\,z_1 z_2 z_3 z_4 = 0   \}
 \, .
$$
The normal direction to $\widehat{X}_5$ on a point is given by the vector
$$
\vec{N} \:=\:( 5z_1^4 +\psi z_2 z_3 z_4 \: ,\:
5z_2^4 + \psi z_1 z_3 z_4 \: ,\:
5z_3^4 + \psi z_1 z_2 z_4 \: ,\:
5z_4^4 + \psi z_1 z_2 z_3 )
$$
For a point\/ $z=(0,0,z_3,z_4) \in \pi_5^{-1}(X_{({\bf g})}^{\circ})$, 
we have\/ $\vec{N}_z =(0,0,5z_3^4,5z_4^4)$.
Hence a local framing for $
T\widehat{X}_{5}\bigr|_{\pi_5^{-1}(X_{({\bf g})}^{\circ})}$ 
is given by tangent vectors
\begin{equation}\label{basis:TX_5}
\left\{ \xi_1 := \frac{\partial}{\partial z_1} \: ,\hspace{6pt}
\xi_2 := \frac{\partial}{\partial z_2} \: ,\hspace{6pt}
\xi_3 := z_4^4 \frac{\partial}{\partial z_3} - z_3^4
\frac{\partial}{\partial z_4}
   \right\}\, .
\end{equation}
$\kG$\/ action on the framing is given by $g_1 (z,\xi_1) = (z,\zeta^n \xi_1)$,
$g_1 (z,\xi_2) = (z,\zeta^{5-n} \xi_2)$, and $g_1 (z,\xi_3) = (z,\xi_3)$.

Now consider the framing
$\{ \xi_1\otimes\overline{\omega}_1 \, ,\,  \xi_2\otimes\overline{\omega}_1 \, ,\,
\xi_3\otimes\overline{\omega}_1 \, ,\, \xi_1\otimes\overline{\omega}_2 \, ,\,
\xi_2\otimes\overline{\omega}_2 \, ,\, \xi_3\otimes\overline{\omega}_2
\}\,$ for the bundle $T\widehat{X}_5 \otimes H^{0,1}(\Sigma)$\/ over
$\pi_5^{-1}(X_{({\bf g})}^{\circ})$.
With respect to this framing, $\lambda_1^{\ast}=g_1$ is given by the matrix
\begin{equation}\label{matrix:6x6}
\left(  \begin{array}{cccccc}
{\scriptstyle \zeta^n (-1-\overline{p}+\overline{q})}&0&0&
{\scriptstyle \zeta^n (-\overline{q}+\overline{s})}&0&0 \\
0&{\scriptstyle \zeta^{5-n} (-1-\overline{p}+\overline{q})}&0&0&
{\scriptstyle \zeta^{5-n} (-\overline{q}+\overline{s}) }&0 \\
0&0& {\scriptstyle -1-\overline{p}+\overline{q}}& 0&0&
{\scriptstyle -\overline{q}+\overline{s}} \\
{\scriptstyle \zeta^n (-\overline{q})}&0&0&
{\scriptstyle \zeta^n (-\overline{s})}&0&0  \\
0&{\scriptstyle \zeta^{5-n} (-\overline{q})}&0&0&
{\scriptstyle \zeta^{5-n} (-\overline{s})}&0  \\
0&0&{\scriptstyle -\overline{q}}&0&0&
{\scriptstyle -\overline{s}}
\end{array}\right)
\end{equation}
where\/ $\zeta =\exp(2\pi i /5)$.  
Matrix~(\ref{matrix:6x6}) can be diagonalized over $\cc$, and we
find that $1$ is an eigenvalue of multiplicity one with corresponding 
eigenvector in Table~\ref{table:eigenvectors}.
\begin{table}[!ht]
\caption{ }
\begin{center}\label{table:eigenvectors}
\begin{tabular}{|c|c|c|c|}
\hline \vspace{-.3cm}
&&&\\
     $n=1$  & $n=2$ & $n=3$ & $n=4$ \\ \vspace{-.3cm}
&&&\\ \hline \vspace{-.3cm} &&&\\ \vspace{-.3cm}
$(0,\zeta^3,0,0,1,0)$  &
$(0,\zeta ,0,0,1,0)$  &
$(\zeta,0,0,1,0,0)$ &
$(\zeta^3,0,0,1,0,0)$ \\
&&&\\ \hline
\end{tabular}
\end{center}
\end{table}
Hence a generator $s$\/ for the restricted obstruction 
bundle\/  $(e^{\ast}TX \otimes H^{0,1}(\Sigma))^{\kG} 
\bigr|_{X_{({\bf g})}^{\circ}}\!\! \rightarrow  X_{({\bf g})}^{\circ} 
$ 
is given by Table~\ref{table:holo-section}.  
Note however that $s$\/ does not define an orbifold section of the restricted 
obstruction bundle.  
\begin{table}[!ht]
\caption{ }
\begin{center}
\label{table:holo-section}
\begin{tabular}{|c|c|}
\hline \vspace{-.3cm}
&\\
\hspace{.3cm} $n=1$ \hspace{.3cm} & 
$\zeta^3(\xi_2\otimes\overline{\omega}_1)+(\xi_2\otimes\overline{\omega}_2)$ \\ \vspace{-.3cm}
&\\ \hline \vspace{-.3cm} &\\ 
\hspace{.3cm} $n=2$ \hspace{.3cm} &
$\zeta(\xi_2\otimes\overline{\omega}_1)+(\xi_2\otimes\overline{\omega}_2)$
 \\ \vspace{-.3cm}
&\\ \hline \vspace{-.3cm} &\\ 
\hspace{.3cm} $n=3$ \hspace{.3cm} &
$\zeta(\xi_1\otimes\overline{\omega}_1)+(\xi_1\otimes\overline{\omega}_2)$ \\ \vspace{-.3cm}
&\\ \hline \vspace{-.3cm} &\\ \vspace{-.3cm}
\hspace{.3cm} $n=4$ \hspace{.3cm} &
$\zeta^3(\xi_1\otimes\overline{\omega}_1)+(\xi_1\otimes\overline{\omega}_2)$ \\ 
&\\   
 \hline
\end{tabular}
\end{center}
\end{table}

 Let $\tilde{p}_3\in \pi^{-1}_5(p_3)$ and
$ \tilde{p}_4 \in \pi^{-1}_5(p_4)$ be any preimages. Then in a neighborhood of
 $\tilde{p}_3$ or $\tilde{p}_4$, 
$(T\widehat{X}_5 \otimes H^{0,1}(\Sigma))^{\kG}$\/ is still generated by the eigenvector 
from  Table~\ref{table:holo-section}.  
Now consider $p_5$ and let $\tilde{p}_5 \in \pi^{-1}_4(p_5)$. Then in a small 
neighborhood of $\tilde{p}_5$ in $\widehat{X}_4$,  
$T\widehat{X}_{4}$ 
is generated by tangent vectors
$$
\left\{ \xi_1' := \frac{\partial}{\partial w_1} \: ,\hspace{6pt}
\xi_2' := \frac{\partial}{\partial w_2} \: ,\hspace{6pt}
\xi_3' := w_5^4 \frac{\partial}{\partial w_3} - w_3^4
\frac{\partial}{\partial w_5}
   \right\}\, .
$$

For any point $w=(0,0,w_3,w_5)$ in such a neighborhood,
$\kG$\/ action on this basis is given by\/ $g_1 (w,\xi_1') = (w,\zeta^n \xi_1')$,
$g_1 (w,\xi_2') = (w,\zeta^{5-n} \xi_2')$, and\/ $g_1 (w,\xi_3') = (w,\xi_3')$.
Hence by an argument completely analogous to the $\widehat{X}_5$ case, a generator
for $(T\widehat{X}_4 \otimes H^{0,1}(\Sigma))^{\kG}$\/ can be obtained  
from  Table~\ref{table:holo-section} by substituting $\xi_j'$ for $\xi_j$.

Now we are in a position to describe the local uniformizing charts for 
$E_{(\bf{g})}$. The following description holds only in the cases $n=1$ or $2$.
The other two cases are completely analogous and left to the reader.   
 
Choose a local chart $(V_x^{{\bf g}} , C({\bf g})  , \pi)$
for $X_{(\bf{g})}$ at a point $x$. 
\begin{itemize}
\item[(i)]
If $x \in X^{\circ}_{(\bf{g})}$, then 
$C({\bf g}) = G_x = \kG \cong \zz_5$, and 
$(V_x^{{\bf g}} \times \cc, \kG ,\tilde{\pi})$ is a uniformizing system for 
$E_{(\bf{g})}$ where $\kG$ acts on 
$V_x^{{\bf g}} \times \cc$ by $g_1(u,v)= (u,\zeta^{5-n}\cdot v)$.
\item[(ii)] If $x=p_3$, then $C({\bf g}) = G_x \cong (\zz_5)^2$. An element of
$C({\bf g})$ is of the form $(\zeta^a,\zeta^b,\zeta^c,1,1)$, where $a+b+c \equiv
0 \mod 5$. Choose generators $g_1$ and $h=(1,\zeta,\zeta^4,1,1)$ of $C({\bf g})$. 
$C({\bf g})$ acts on $V_x^{{\bf g}} \times \cc$ as follows:
$g_1(u,v)= (u,\zeta^{5-n}\cdot v)$, and $h(u,v)= (\zeta^4\cdot u,\zeta\cdot v)$.
\item[(iii)] If $x=p_4$, then $C({\bf g})=\{(\zeta^a, \zeta^b,1,\zeta^c,1)
\:|\: a+b+c \equiv
0 \mod 5 \}$.   Choose $h=(1,\zeta,1,\zeta^4,1)$.  Then   
the $C({\bf g})$ action on $V_x^{{\bf g}} \times \cc$\/ is 
given by the same formulas as for $p_3$.
\item[(iv)]  If $x=p_5$, then $C({\bf g})=\{(\zeta^a, \zeta^b,1,1,\zeta^c )
\:|\: a+b+c \equiv
0 \mod 5 \}$.   Choose $h=(1,\zeta,1,1,\zeta^4)$.  
Then   
the $C({\bf g})$ action on $V_x^{{\bf g}} \times \cc$\/ is 
given by the same formulas as for $p_3$.
\end{itemize}
We conclude that there is a smooth $\kG$ action on $E_{(\bf{g})}$ that preserves
the fiber.

$E_{(\bf{g})}$ is a line bundle. So the $3$-point function (\ref{eq:3point}) is
nonzero only if $\eta_i \in H^{0,0} (X_{(g_i)})$, $ 1\leq i\leq 3$. If this is the case
then  
\begin{equation}\label{eq:eta-product}
\langle \eta_1, \eta_2, \eta_3 \rangle_{orb} \:=\:
\eta_1\eta_2\eta_3\cdot 
\langle c_1(E_{(\bf{g})}),[X_{(\bf{g})}]\rangle 
\:=\: \eta_1\eta_2\eta_3 \int_{X_{(\bf{g})}}^{orb} c_1(A)
\end{equation} 
where $A$\/
is any orbifold connection on $E_{(\bf{g})}$.

Now consider the associated orbifold 
principal $S^1$ bundle $P_{({\bf g})}$ such that 
$E_{(\bf{g})} = P_{({\bf g})}\times_{S^1} \cc$.  
Recall that there is a global action of $\kG\cong \zz_5$\/ 
on each fiber $F=S^1$ of $P_{({\bf g})}$.  The quotient 
$P_{({\bf g})}/\kG$\/ is again an orbifold  
principal bundle over the orbifold $X_{({\bf g})}$.  
Let\/ 
$\pi_{\kG} :P_{({\bf g})} \rightarrow  P_{({\bf g})}/\kG$\/ be the quotient map, which 
extends to an orbifold bundle map.  Choose
an orbifold connection $A$\/ that is the pullback $\pi_{\kG}^{\ast}(A')$, where $A'$ 
is an orbifold connection on the associated bundle $E'_{({\bf g})}=
(P_{({\bf g})}/\kG)\times_{S^1}\cc$ over $X_{({\bf g})}$. 
Note that $\pi_{\kG}$ on each fiber is given by $z\mapsto z^5$.  
The Lie algebra of $F$\/ can be identified with $\rr$.  
Hence the induced map on the Lie algebra 
$(\pi_{\kG})_{\ast}:\rr\rightarrow \rr$ is just mutiplication by 5.  

Let $\Omega$ and $\Omega'$ be curvature 2-forms for $A$\/ and $A'$.  
By Proposition~6.1 of \cite{[KN]}, we must have\/  
$\pi_{\kG}^{\ast}(\Omega') = 5 \cdot \Omega$.  
Hence by Chern-Weil Theory, we have 
\begin{equation}\label{eq:A=A'}
\int_{X_{(\bf{g})}}^{orb} c_1(A) = \frac{1}{5}
\int_{X_{(\bf{g})}}^{orb} c_1(A') \, .
\end{equation}

Since the action of $\kG$ in any uniformizing system of $E'_{(\bf{g})}$ is trivial,
$E'_{(\bf{g})}$ induces an orbifold bundle $E''_{(\bf{g})}$ over the reduced orbicurve
$X_{(\bf{g})}'$ which has an induced connection $A''$.  The connections $A'$
and $A''$ may be represented by the same $1$-form on $V$\/ for any pair of corresponding
uniformizing systems $(V \times \cc, G', \tilde{\pi}'_1)$ 
and $(V \times \cc, G'/\kG, \tilde{\pi}''_1)$
of $E'_{(\bf{g})}$ and $E''_{(\bf{g})}$ respectively. 
By Chern-Weil theory, $c_1(A')$ and $c_1(A'')$
may therefore be represented by the same $2$-form on $V_1$. Hence by 
(\ref{eq:orb-integral}) 
\begin{equation}\label{eq:A'=A"}
\int_{X_{(\bf{g})}}^{orb} c_1(A')\:=\:\frac{1}{5}\int_{X_{(\bf{g})}'}^{orb} c_1(A'')
\:=\: \frac{1}{5} \langle c_1(E''_{({\bf g})}), [X_{(\bf{g})}'] \rangle .
\end{equation}
Finally we know that
(cf.$\;$\cite{[CR1]}, equation~(4.2.5)\/)
\begin{equation}\label{eq:desingularization}
\langle c_1(E''_{({\bf g})}), [X_{({\bf g})}']\rangle \:=\:
\langle c_1( | E''_{({\bf g})} |), [\Sigma_{0}]\rangle \: +\:
\sum_{j=3}^{5}\frac{\mu_j}{5} \; ,
\end{equation}
where $| E''_{({\bf g})}|$ denotes the desingularization of\/ $E''_{({\bf g})}$ 
over a smooth genus $0$ curve $\Sigma_0$,  and\/
$0\leq \mu_j < 5$\/ are the
local invariants of $E''_{({\bf g})}$ at the singular points $p_j$
of $X_{({\bf g})}'$.

To pin down the local invariants $\mu_j$, we need to describe the
uniformizing system for $E''_{({\bf g})}$.  
Choose a local chart $(V_x^{{\bf g}} , C({\bf g})/\kG  , \pi'')$
for $X_{(\bf{g})}'$ at a point $x$. In the $n=1$ or $2$ cases,
\begin{itemize}
\item[(i)]  If\/ $x \in X^{'\circ}_{(\bf{g})}$, then\/ 
$C({\bf g}) / \kG$ is the trivial group.   
\item[(ii)]  If\/ $x=p_3$, then\/ $C({\bf g}) / \kG\cong \zz_5$ and $\zz_5$ acts on
 $V_x^{{\bf g}} \times \cc$ by $h(u,v) = (\zeta^4\cdot u,v)$. 
\end{itemize}
From (ii) we immediately conclude that all $\mu_j=0$.

To compute $c_1(|E''_{({\bf g})}|)$, note that $s^5$ 
is a holomorphic orbifold section of $ E''_{({\bf g})}
\bigr|_{X_{({\bf g})}^{'\circ}}\!\!\rightarrow 
X_{({\bf g})}^{'\circ}$.  
This section can be extended to a global holomorphic section of 
the desingularized bundle $|E''_{({\bf g})}|$ as follows.

Note that we can express our framing (\ref{basis:TX_5}) 
in terms of $w_i$ coordinates.  
On $V_4$ we have
\begin{align}\label{eq:xi_1}
\xi_1 = \frac{\partial}{\partial z_1} & = \sum_{j=1}^{5}
\frac{\partial w_j}{\partial z_1} \cdot\frac{\partial}{\partial w_j}
=\sum_{j=1}^{5} \frac{\partial}{\partial z_1} \!\left(
\frac{z_j}{z_4} \right)\cdot \frac{\partial}{\partial w_j}
\\[5pt]
 & = \frac{1}{z_4} \cdot \frac{\partial}{\partial w_1}
  =
\frac{w_5}{w_4} \cdot \frac{\partial}{\partial w_1} =
w_5 \, \frac{\partial}{\partial w_1} \,  .\notag
\end{align}
Similarly, we have
\begin{equation}\label{eq:xi_2}
 \xi_2 \: = \frac{\partial}{\partial z_2} =
\sum_{j=1}^{5}
\frac{\partial w_j}{\partial z_2} \cdot\frac{\partial}{\partial w_j}
= \frac{1}{z_4} \cdot \frac{\partial}{\partial w_2}
= w_5 \frac{\partial}{\partial w_2} \,  .
\end{equation}
One can express $\xi_3$ similarly, but we will not need it.   

Now note that $s^5$ defines a global holomorphic section of\/ 
$| E''_{({\bf g})}|$, where
\begin{equation*}
s^5\bigr|_{U_4 \cap \Sigma_{0}} = \:(\zeta w_5)^5 \,\frac{\partial}{\partial
w_2}\otimes\overline{\omega}_1  \:+\:
(w_5)^5 \,\frac{\partial}{\partial w_2}\otimes\overline{\omega}_2 \: .
\end{equation*}
We also note that\/  
$(w_5)^5$ defines the local coordinate of the desingularized curve 
$\Sigma_0$ centered at\/ $p_5\,$.  
Hence the section\/ $s^5$ has a unique zero of order one at the point 
$\,p_5=[x_1:x_2:x_3:x_4:x_5]=[0:0:-\zeta:1:0]\,$
in $\Sigma_{0}$.

Since a generic holomorphic section of\/ $| E''_{({\bf g})}|$\/  will also 
have a single zero, 
we conclude that
$$
\langle c_1( | E''_{({\bf g})} |), [\Sigma_{0}] \rangle \:=\: 1 
\, .
$$
It follows from (\ref{eq:A=A'}), (\ref{eq:A'=A"}) and (\ref{eq:desingularization}) that 
$$
\left\langle c_1(E_{({\bf g})}), [X_{({\bf g})}]\right\rangle \:=\:\frac{1}{25} \,  .
$$
Finally, (\ref{eq:eta-product}) becomes 
\begin{equation}\label{final answer}
\langle \eta_1, \eta_2, \eta_3 \rangle_{orb} \:=\:
\frac{1}{25} \, \eta_1\eta_2\eta_3 \: .
\end{equation}
We remark that (\ref{final answer}) holds true when\/ $n=3$ or $4$, 
as well as for the other choices 
of $({\bf g})$, as long as $g_2 \neq (g_1)^4$.

\subsubsection{The Case\/ $g_2=(g_1)^2$}
\begin{figure}[!ht]
\begin{center}
\includegraphics[scale=.6]{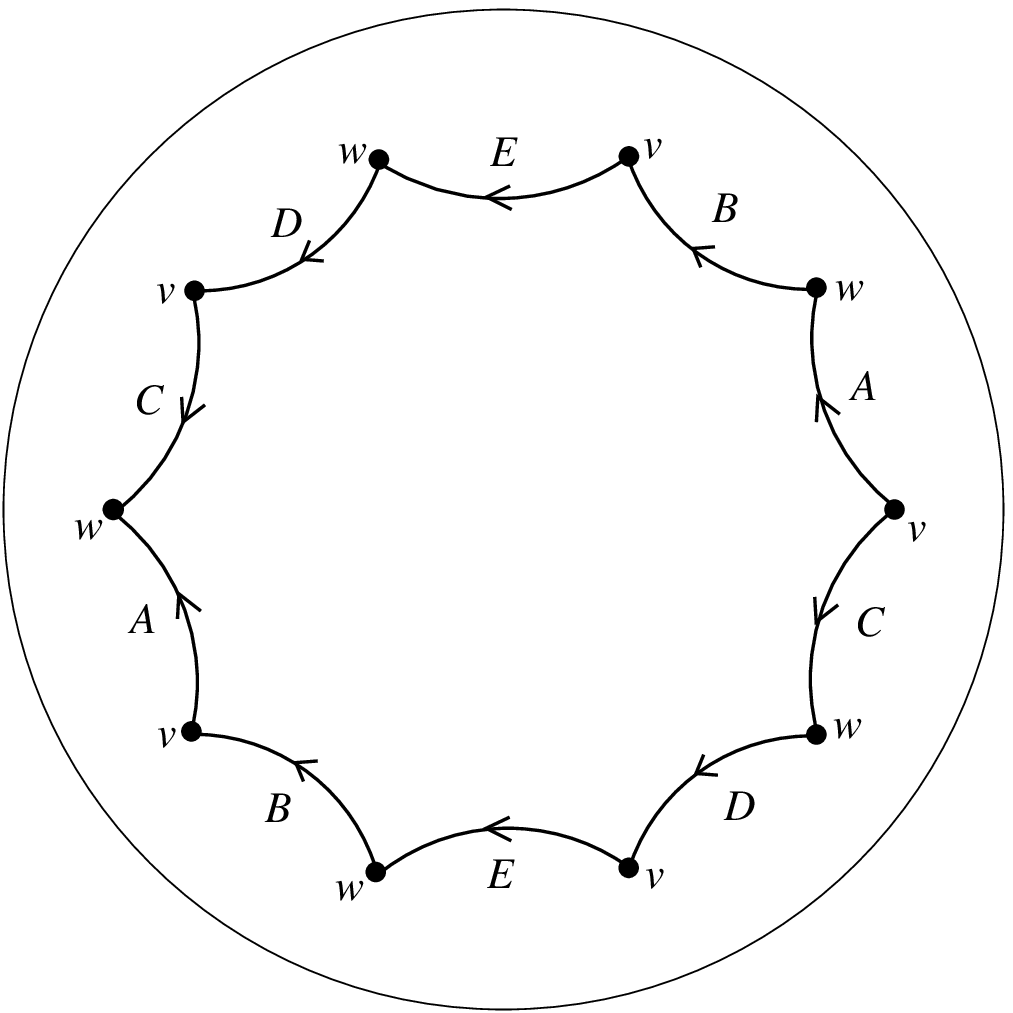}
\end{center}
\caption{ }
\label{fig:tengon2}
\end{figure}

\ For this choice of $\rho$, the quotient\/ $\Sigma = \hh^2 /{\rm Ker}(\rho)$\/ is
the genus two surface given by the decagon in Figure~\ref{fig:tengon2} whose sides are
identified according to Table~\ref{table:tengon:k=2}.

\begin{table}[!ht]
\caption{ }
\begin{center}\label{table:tengon:k=2}
\begin{tabular}{|c|c|c|c|c|}
\hline \vspace{-.3cm}
&&&&\\
     $A$  & $B$ & $C$ & $D$ & $E$\\ \vspace{-.3cm}
&&&&\\ \hline \vspace{-.3cm} &&&&\\ \vspace{-.3cm}
\hspace{10pt}$\lambda_1^{-2}\lambda_2$\hspace{10pt}  &
\hspace{10pt}$\lambda_2\lambda_1^{-2}$\hspace{10pt} &
$\lambda_1^3(\lambda_1^{-2}\lambda_2)\lambda_1^{-3}$ &
$\lambda_1(\lambda_1^{-2}\lambda_2)\lambda_1^{-1} $
& $\lambda_1^4(\lambda_1^{-2}\lambda_2)\lambda_1^{-4} $ \\
&&&&\\ \hline
\end{tabular}
\end{center}
\end{table}

After cutting and pasting along the geodesics as in
the previous case, we can find a set of
generators of $\pi_1(\Sigma)$ based at the point\/ $v\,$:
$$
\alpha = AB, \quad \beta=EDCB, \quad \gamma = ED, \quad \delta=
CE^{-1}
$$
which descends to a canonical symplectic basis for $H_1(\Sigma;\zz)$.
The rotation $\lambda_1$ induces the following map on
$\pi_1(\Sigma)$:
$$
\begin{array}{l}
\alpha \:\longmapsto\: D^{-1}C^{-1} = \gamma^{-1} \delta^{-1} \\[7pt]
\beta \:\longmapsto\: ABEC^{-1} = \alpha \,\delta^{-1}  \\[7pt]
\gamma \:\longmapsto\: AB = \alpha \\[7pt]
\delta \:\longmapsto\: EA^{-1}= \delta^{-1}\gamma^{-1}
\beta\alpha^{-1}
\end{array}
$$
Hence the automorphism $\lambda_1^{\ast}:H^1(\Sigma)\rightarrow
H^1(\Sigma)$ can be expressed in the matrix
\begin{equation}\label{matrix:H^1:k=2}
\left(\begin{array}{rrrr}0&1&1&-1 \\
0&0&0&1 \\
-1&0&0&-1 \\
-1&-1&0&-1
\end{array}\right)^T
\;=\;\left(\begin{array}{rrrr}
0&0&-1&-1 \\
1&0&0&-1 \\
1&0&0&0 \\
-1&1&-1&-1
\end{array}\right)
\end{equation}
with respect to the dual basis $\left\{\check{\alpha},
\check{\beta}, \check{\gamma}, \check{\delta} \right\}$.

Let $\,\{ \omega_1 , \omega_2 \} = \{ \check{\alpha} + p\check{\beta} + q\check{\delta} \, ,\, \check{\gamma}+q\check{\beta}+s\check{\delta} \}\,$ be a basis for $H^{1,0}(\Sigma)$ as before.
From matrix (\ref{matrix:H^1:k=2}), we calculate that
$$
\begin{array}{l}
\lambda_1^{\ast} \omega_1 \:=\: (-q) \check{\alpha}
+(1-q)\check{\beta} + \check{\gamma} +(-1+p-q)\check{\delta} \:=\:
 (-q)\omega_1 + \omega_2 \: ,\\[7pt]
\lambda_1^{\ast}\omega_2 \:=\:  (-1-s) \check{\alpha}
+(-s)\check{\beta}  +(-1+q-s)\check{\delta} \:=\: (-1-s)\omega_1 \: .
\end{array}
$$
Comparing coefficients, we obtain the system:
\begin{equation}\label{eq:pqs:k=2}
\left\{ \begin{array}{l}
1-q = -qp+q \\[7pt]
-1+p-q= -q^2+s\\[7pt]
-s= p(-1-s) \\[7pt]
-1+q-s=q(-1-s)
\end{array}\right.
\end{equation}
Just as in the previous case,
there are four solutions to system (\ref{eq:pqs:k=2}), but only
one of them satisfies the condition that Im$P$ is positive definite, namely,

\begin{equation}\label{pqs:k=2}
\left\{\begin{array}{l}
p\:=\:  \frac{5-\sqrt{5}}{4}\,+\, i \sqrt{\frac{5+\sqrt{5}}{8}} \:= \:
\sqrt{\frac{5-\sqrt{5}}{2}} \exp(3\pi i /10) \, ,
\\[7pt]
q\:=\: \frac{1}{2}\left( 1+i\sqrt{5-2\sqrt{5}} \right)  \:=\:
2\sin(\pi /10) \cdot \exp(\pi i /5)\, , \\[7pt]
s\:=\:  \frac{-5+\sqrt{5}}{4}\,+\, i \sqrt{\frac{5+\sqrt{5}}{8}} \:=\:
\sqrt{\frac{5-\sqrt{5}}{2}} \exp(7\pi i /10)\,  .
\end{array}
\right.
\end{equation}

\subsubsection*{Representative Subcase}
Without loss of generality, we choose\/
$g_1=(\zeta^n ,\zeta^{5-n},1,1,1)$, where $n=1,2,3,4$.
Once again,
consider the framing
\begin{equation}\label{basis:e*TX otimes H^0,1}
\{ \xi_1\otimes\overline{\omega}_1 \, ,\,  \xi_2\otimes\overline{\omega}_1 \, ,\,
\xi_3\otimes\overline{\omega}_1 \, ,\, \xi_1\otimes\overline{\omega}_2 \, ,\,
\xi_2\otimes\overline{\omega}_2 \, ,\, \xi_3\otimes\overline{\omega}_2\}\,
\end{equation}
for the bundle $\,T\widehat{X}_5 \otimes H^{0,1}(\Sigma)$\/ over
$\,\pi_5^{-1}(X^{\circ}_{({\bf g})})$.
With respect to this framing, $\,\lambda_1^{\ast}=g_1$ is given by the matrix
\begin{equation}\label{matrix:6x6:k=2}
\left(  \begin{array}{cccccc}
{\scriptstyle -\zeta^n\,\overline{q}}&0&0&
{\scriptstyle \zeta^n(-1-\overline{s})}&0&0 \\
0&{\scriptstyle -\zeta^{5-n}\,\overline{q}}&0&0&
{\scriptstyle \zeta^{5-n}(-1-\overline{s}) }&0 \\
0&0& {\scriptstyle -\overline{q}}& 0&0&
{\scriptstyle -1-\overline{s}} \\
{\scriptstyle \zeta^n}&0&0&
0&0&0  \\
0&{\scriptstyle \zeta^{5-n}}&0&0&
0&0  \\
0&0&1&0&0&0
\end{array}\right)
\end{equation}
where\/ $\zeta = \exp(2\pi i /5)$.
Matrix (\ref{matrix:6x6:k=2}) is diagonalizable over $\cc$, and
one finds that $1$ is an eigenvalue of multiplicity one with corresponding
eigenvector in Table~\ref{table:eigenvectors:k=2}.
\begin{table}[!ht]
\caption{ }
\begin{center}\label{table:eigenvectors:k=2}
\begin{tabular}{|c|c|c|c|}
\hline \vspace{-.3cm}
&&&\\
     $n=1$  & $n=2$ & $n=3$ & $n=4$ \\ \vspace{-.3cm}
&&&\\ \hline \vspace{-.3cm} &&&\\ \vspace{-.3cm}
$(0,\zeta,0,0,1,0)$  &
$(\zeta^3,0,0,1,0,0)$ &
$(0,\zeta^3,0,0,1,0)$ &
$(\zeta ,0,0,1,0,0)$ \\
&&&\\ \hline
\end{tabular}
\end{center}
\end{table}
The rest of the computation goes exactly like in the
previous case almost verbatim.

\subsubsection{The Case\/ $g_2=(g_1)^3$}
\begin{figure}[!ht]
\begin{center}
\includegraphics[scale=.6]{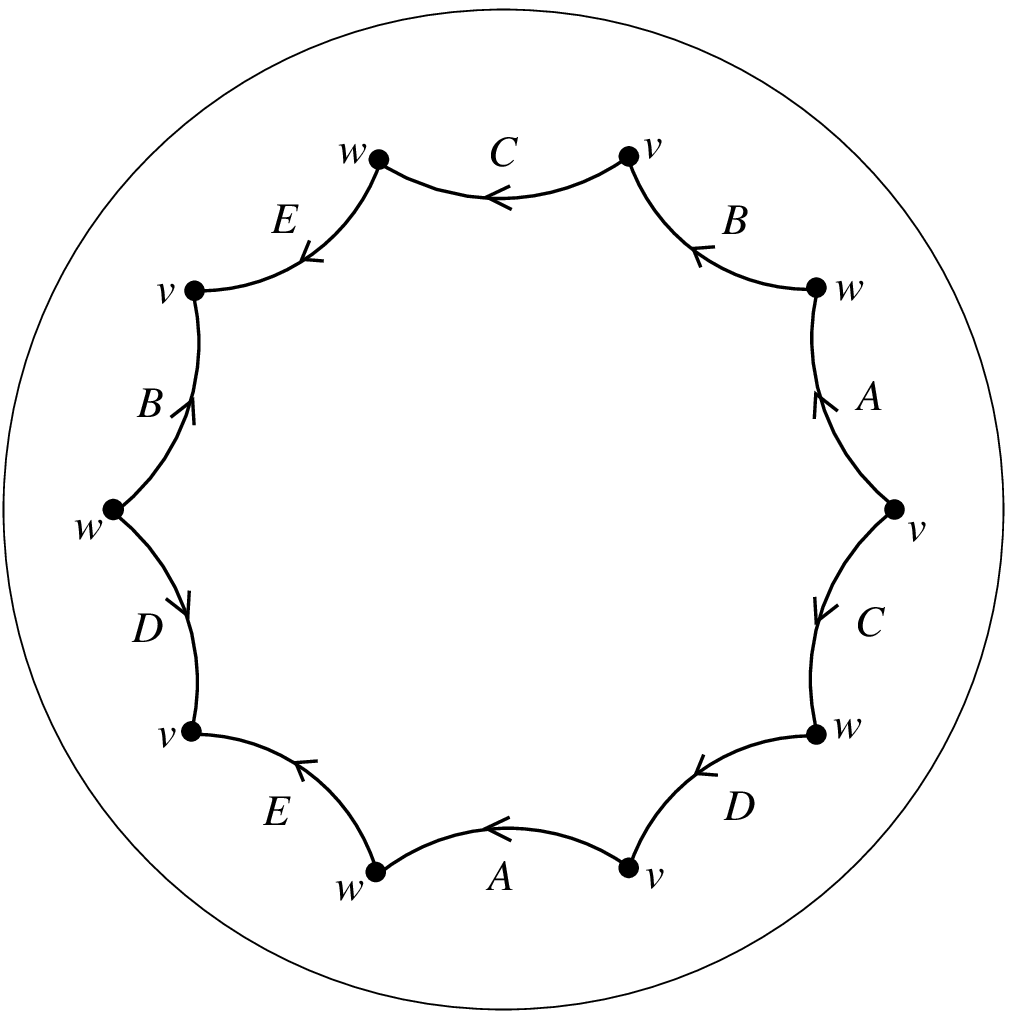}
\end{center}
\caption{ }
\label{fig:tengon3}
\end{figure}

\ For this choice of $\rho$, the quotient\/ $\Sigma = \hh^2 /{\rm Ker}(\rho)$\/ is
the genus two surface given by the decagon in Figure~\ref{fig:tengon3} whose sides are
identified according to Table~\ref{table:tengon:k=3}.

\begin{table}[!ht]
\caption{ }
\begin{center}\label{table:tengon:k=3}
\begin{tabular}{|c|c|c|c|c|}
\hline \vspace{-.3cm}
&&&&\\
     $A$  & $B$ & $C$ & $D$ & $E$\\ \vspace{-.3cm}
&&&&\\ \hline \vspace{-.3cm} &&&&\\ \vspace{-.3cm}
\hspace{10pt}$\lambda_1^{-3}\lambda_2$\hspace{10pt}  &
\hspace{10pt}$\lambda_2\lambda_1^{-3}$\hspace{10pt} &
$\lambda_1^4(\lambda_1^{-3}\lambda_2)\lambda_1^{-4}$ &
$\lambda_1(\lambda_1^{-3}\lambda_2)\lambda_1^{-1} $
& $\lambda_1^2(\lambda_1^{-3}\lambda_2)\lambda_1^{-2} $ \\
&&&&\\ \hline
\end{tabular}
\end{center}
\end{table}

After cutting and pasting along the geodesics as in
previous cases, we can find a set of
generators of $\pi_1(\Sigma)$ based at the point\/ $v\,$:
$$
\alpha = E^{-1}C^{-1}, \quad \beta=AB, \quad \gamma = AD, \quad \delta=
E^{-1}A^{-1}
$$
which descends to a canonical symplectic basis for $H_1(\Sigma;\zz)$.
The rotation $\lambda_1$ induces the following map on
$\pi_1(\Sigma)$:
$$
\begin{array}{l}
\alpha \:\longmapsto\: B^{-1}A^{-1} = \beta^{-1}\\[7pt]
\beta \:\longmapsto\: D^{-1}C^{-1} = \gamma^{-1}\delta^{-1}\alpha  \\[7pt]
\gamma \:\longmapsto\: D^{-1}E = \gamma^{-1}\delta^{-1} \\[7pt]
\delta \:\longmapsto\: B^{-1}D = \beta^{-1}\gamma
\end{array}
$$
Hence the automorphism $\lambda_1^{\ast}:H^1(\Sigma)\rightarrow
H^1(\Sigma)$ can be expressed in the matrix
\begin{equation}\label{matrix:H^1:k=3}
\left(\begin{array}{rrrr}0&1&0&0 \\
-1&0&0&-1 \\
0&-1&-1&1 \\
0&-1&-1&0
\end{array}\right)^T
\;=\;\left(\begin{array}{rrrr}
0&-1&0&0 \\
1&0&-1&-1 \\
0&0&-1&-1 \\
0&-1&1&0
\end{array}\right)
\end{equation}
with respect to the dual basis $\left\{\check{\alpha},
\check{\beta}, \check{\gamma}, \check{\delta} \right\}$.

Let $\,\{ \omega_1 , \omega_2 \} = \{ \check{\alpha} + p\check{\beta} + q\check{\delta} \, ,\, \check{\gamma}+q\check{\beta}+s\check{\delta} \}\,$ be a basis for $H^{1,0}(\Sigma)$ as before.
From
matrix (\ref{matrix:H^1:k=3}), we calculate that
$$
\begin{array}{l}
\lambda_1^{\ast} \omega_1 \:=\: (-p) \check{\alpha}
+(1-q)\check{\beta} + (-q)\check{\gamma} +(-p)\check{\delta} \:=\:
(-p)\omega_1 + (-q)\omega_2  \: ,
\\[7pt]
\lambda_1^{\ast}\omega_2 \:=\:  (-q) \check{\alpha}
+(-1-s)\check{\beta} + (-1-s)\check{\gamma} +(1-q)\check{\delta} \:=\:
(-q)\omega_1 +(-1-s)\omega_2 \: .
\end{array}
$$
Comparing coefficients, we obtain the system:
\begin{equation}\label{eq:pqs:k=3}
\left\{ \begin{array}{l}
1-q = -p^2-q^2 \\[7pt]
-p= -pq-qs\\[7pt]
-1-s= -pq+q(-1-s) \\[7pt]
1-q=-q^2 + s(-1-s)
\end{array}\right.
\end{equation}
which has a unique solution satisfying the positive definiteness condition, namely,
\begin{equation}\label{pqs:k=3}
\left\{\begin{array}{l}
p\:=\: \frac{1}{10}\, i\,( 5+3\sqrt{5} )\sqrt{5-2\sqrt{5}} \: ,
\\[7pt]
q\:=\: \frac{1}{10}\left(5-i\sqrt{5(5-2\sqrt{5})}\right) \:=\:
\sqrt{\frac{5-\sqrt{5}}{10}} \, \exp(-\pi i /10) \, , \\[9pt]
s\:=\: \frac{1}{10}\left( -5 + i\sqrt{5(5+2\sqrt{5})} \right) \:=\:
\sqrt{\frac{5+\sqrt{5}}{10}} \, \exp(7 \pi i /10)\, .
\end{array}
\right.
\end{equation}

\subsubsection*{Representative Subcase}
Without loss of generality, we choose\/
$g_1=(\zeta^n,\zeta^{5-n},1,1,1)$, where\/ $n=1,2,3,4$.
With respect to the framing (\ref{basis:e*TX otimes H^0,1}) for\/
 $T\widehat{X}_5 \otimes H^{0,1}(\Sigma)$\/ over 
$\pi_5^{-1}(X^{\circ}_{({\bf g})})$, 
$\lambda_1^{\ast}=g_1$ is given by the matrix
\begin{equation}\label{matrix:6x6:k=3}
\left(  \begin{array}{cccccc}
{\scriptstyle -\zeta^n\,\overline{p}}&0&0&
{\scriptstyle -\zeta^n\,\overline{q}}&0&0 \\
0&{\scriptstyle -\zeta^{5-n}\,\overline{p}}&0&0&
{\scriptstyle -\zeta^{5-n}\,\overline{q} }&0 \\
0&0& {\scriptstyle -\overline{p}}& 0&0&
{\scriptstyle -\overline{q}} \\
{\scriptstyle -\zeta^n\,\overline{q}}&0&0&
{\scriptstyle \zeta^n (-1-\overline{s})}&0&0  \\
0&{\scriptstyle -\zeta^{5-n}\,\overline{q}}&0&0&
{\scriptstyle \zeta^{5-n} (-1-\overline{s})}&0  \\
0&0&{\scriptstyle -\overline{q}}&0&0&{\scriptstyle -1-\overline{s}}
\end{array}\right)
\end{equation}
Matrix (\ref{matrix:6x6:k=3}) is diagonalizable, and we find that $1$ is
an eigenvalue of multiplicity one with corresponding
eigenvector in Table~\ref{table:eigenvectors:k=3}.
\begin{table}[!ht]
\caption{ }
\begin{center}
\label{table:eigenvectors:k=3}
\begin{tabular}{|c|c|}
\hline \vspace{-.3cm}
&\\
\hspace{.3cm} $n=1$ \hspace{.3cm} & 
$(0,2\,{\rm Re}(\zeta^2),0,0,1,0)$ \\ \vspace{-.3cm}
&\\ \hline \vspace{-.3cm} &\\ 
\hspace{.3cm} $n=2$ \hspace{.3cm} &
$(0,2\,{\rm Re}(\zeta)  ,0,0,1,0)$ \\ \vspace{-.3cm}
&\\ \hline \vspace{-.3cm} &\\ 
\hspace{.3cm} $n=3$ \hspace{.3cm} &
$(2\,{\rm Re}(\zeta)  ,0,0,1,0,0)$ \\ \vspace{-.3cm}
&\\ \hline \vspace{-.3cm} &\\ \vspace{-.3cm}
\hspace{.3cm} $n=4$ \hspace{.3cm} &
$(2\,{\rm Re}(\zeta^2),0,0,1,0,0)$ \\ 
&\\   
 \hline
\end{tabular}
\end{center}
\end{table}

The rest of the computation goes exactly like in the
previous two cases.

\subsubsection{The Case\/ $g_2=(g_1)^4$}
As we remarked earlier in this case $\Sigma =\sss^2$. Hence $H^{0,1}(\Sigma)$ is trivial
and $E_{(\bf{g})}$ is a rank zero bundle. So $e_A(E_{(\bf{g})}) =1$. Thus 
the $3$-point function (\ref{eq:3point}) 
$$\langle \eta_1 , \eta_2 , \eta_3 \rangle_{orb} \:= \:
\int_{X_{({\bf g})}}^{orb} e_1^{\ast} \eta_1 \wedge e_2^{\ast}\eta_2
\wedge e_3^{\ast}\eta_3 .$$
For the integral to be nonzero exactly one of the $\eta_j$ must be a $(1,1)$-form 
and the other two must be $0$-forms.

Let $H$ denote the hyperplane class in
$H_6(\cpk/\zz_5^3)$.  Let $\eta \in H^{1,1}(X)$ denote the generator whose
Poincar\'e dual is $[H\cap X]\in H_4(X)$. Then $\eta'=e_3^{\ast}\eta$ generates
$H^{1,1}(X_{({\bf g})})$. Moreover, the maps $e_1$ and $e_2$ are isomorphisms
in this case. Hence we only need to compute $ \int_{X_{({\bf g})}}^{orb} \eta'\/.$

First write $e_3$ as the composition\/ $X_{({\bf g})}\! \stackrel{q}{\longrightarrow}\!
X_{({\bf g})}' \stackrel{j}{\hookrightarrow} X$, where $q$ is the reduction map.  
We have 
$$ \int_{X_{({\bf g})}}^{orb} \! e_3^{\ast}\eta \;=\; \frac{1}{5} 
\int_{X_{({\bf g})}'}^{orb} j^{\ast}\eta 
\;=\; \frac{1}{5} \int_{j(X_{({\bf g})}')}^{orb} \!\eta 
\;=\; \frac{1}{5} \int_{j(X_{({\bf g})}')} \!\eta \: .
$$
Note that the last equality follows from the fact that\/ 
$X_{({\bf g})}'$\/ is reduced.  

Since $X$\/ is a normal complex algebraic variety, there is a
well-defined notion of Poincar\'e duality (cf.$\;$\cite{[GM]}).
Hence the last integral is given by the intersection number of 
$[j(X_{({\bf g})}')]$ and $PD(\eta)$; it is easily 
checked that this is $1$. 
Hence $ \int_{X_{({\bf g})}}^{orb} \eta'\/=\frac{1}{5}$.

\section{Ordinary Cup Product on Mirror Quintic}\label{sec:ord-cup}

First we determine the Hodge numbers of $X$.  We start with
Lefschetz Hyperplane Theorem in \cite{[BC]} (Proposition 10.8).
\begin{theorem} Let $X$ be a nondegenerate
 ample hypersurface of an $n$-dimensional complete simplicial toric variety
 $Y$. Then the natural map induced by inclusion
 $j^{*}: H^i (Y) \to H^i (X)$, is an isomorphism for\/
 $i<n-1$ and an injection for\/ $i=n-1$.
\end{theorem}

Applying the above theorem to our mirror quintic $X$ inside
$Y=\cpk /(\zz_5)^3$, we easily conclude that $h^{1,1}(X)=1\,$
and\/ $h^{1,0}(X)=h^{2,0}(X)=0$. Complex conjugation gives
$h^{0,1}(X)=h^{0,2}(X)=0$.  By Kodaira-Serre duality
(cf.$\;$\cite{[CR1]}, Proposition~3.3.2), it follows that
$h^{3,2}(X)=h^{3,1}(X)=0\,$ and\/ $h^{2,2}(X)=1$. By Theorem~4 in
\cite{[Po]}, we can calculate that $h^{2,1}(X)=1$.

By Poincar\'e duality, we may compute the
cup product by computing the intersection of dual cycles
(cf.$\;$\cite{[GM]}). Let $H$ denote the hyperplane class in
$H_6(Y)$.  Let $\eta \in H^{1,1}(X)$ denote the generator whose
Poincar\'e dual is $[H\cap X]\in H_4(X)$.  To compute $\eta \cup
\eta \cup \eta$, we look at the intersection of $X$\/ and three
hyperplanes in general position.  In particular, we look at
hyperplanes\/ $x_1 =0$, $x_2=0$, and\/ $x_3=0$. These intersect
$X$\/ transversally at the unique point $[0:0:0: -\zeta : 1]\,$ in
$Y$. Thus $\eta\cup\eta\cup\eta =PD([X])$.

\subsection{Middle Cohomology}
The cup product for the middle cohomology of an ample  
hypersurface is described for the smooth case in \cite{[CG]} and for 
the quasi-smooth case in \cite{[Ma]}.

First we state some general facts. Let $X$ be a nondegenerate
ample hypersurface in a $d$-dimensional projective simplicial
toric variety $Y$. Let $S = \cc[x_1,\ldots\hspace{-1pt},x_n]$ be the homogeneous
coordinate ring of\/ $Y$ which is graded by the Chow group $A_{d-1}(Y)$.
Suppose $X$ is defined by a degree $\beta$ homogeneous polynomial $f \in S_\beta$.  
Denote by $D_i$ the effective
divisor $\{x_i=0\}$ and let $D=\sum b_iD_i$ be an ample divisor such that 
$[D]=\beta$. Set $\beta_0 = \deg(x_1\ldots x_n)=[\sum_{i=1}^n D_i]$. 

Let $F_j = x_j\frac{\partial f}{\partial x_j}$. 
Given $f \in S_\beta $ we get the ideal quotient 
(cf.$\;$\cite{[CLO]}, p.$\,$191) $$ J_1(f) = \lan F_1,\ldots\hspace{-1pt},F_n \ran 
 : (x_1\ldots x_n) $$
 and the ring $R_1(f) = S/J_1(f)$ graded by $A_{d-1}(Y)$.

Fix an integer basis $m_1,\ldots\hspace{-1pt},m_d$ for the lattice $M$. 
Moreover suppose $e_1,\ldots\hspace{-1pt},e_n$ are the primitive 
integral generators of the $1$-dimensional cones of the fan of $Y$. Then given a 
subset $I = \{i_1,\ldots\hspace{-1pt},i_d \} \subset \{1,\ldots\hspace{-1pt},n\}$,
we define\/ $$\det(e_I) = \det( \lan m_j, e_{i_k} \ran_{1\leq j\leq d,\,
i_k \in I} ),$$
$dx_I = dx_{i_1}\wedge \hspace{-1pt}\ldots \wedge dx_{i_d}\,$, and $\hat{x}_I =
\prod_{i \not\in I} x_i \,$.
Define the $n$-form $\Omega$ by the formula 
$$ \Omega = \sum_{|I|=d} \det(e_I)\hat{x}_I dx_I .$$

For $A \in S_{(a+1)\beta - \beta_0}$ consider the rational $d$-form
$\omega_A = A \Omega / f^{a+1}$ which gives a class in $H^d(Y-X)$ and
by the residue map $$ {\rm Res} : H^d(Y-X) \to H^{d-1}(X) $$ we get
${\rm Res}(\omega_A) \in H^{d-1}(X)$. 
 Denote the $(b,a)$ Hodge 
component of ${\rm Res}(\omega_A)$ by ${\rm Res}(\omega_A)^{b,a}$,
where $a+b = d-1$.
Then the map 
\begin{equation}\label{map:residue}
{\rm Res}(\omega_{-})^{b,a} : R_1(f)_{(a+1)\beta -\beta_0}
\to  H^{b,a}(X)
\end{equation}
is injective by Theorem~4.4 of \cite{[Ma]}. 
By Theorem~11.8 
of \cite{[BC]}, the vector spaces $R_1(f)_{(a+1)\beta -\beta_0}$ and 
$PH^{b,a}(X) := (H^{d-1}(X)/j^{\ast} H^{d-1}(Y))^{b,a}$ 
have the same dimension. If $d$\/ is even
then $H^{d-1}(Y)= 0$, and hence\/ $PH^{b,a}(X) = H^{b,a}(X)$. Thus for even $d$,
which is true for our mirror quintic case, the map (\ref{map:residue})
is an isomorphism. 

Now let $I=\{i_0,\ldots\hspace{-1pt},i_d\} \subset \{1,\ldots\hspace{-1pt},n \}$.
Let\/ $J=\det(\frac{\partial F_j}{\partial x_i})_{i,j \in I}/
(c_I^\beta)^2 \hat{x}_I$, where $(c_I^\beta)$ is the determinant of the 
$(d+1) \times (d+1)$ matrix obtained from $( \lan m_j, e_{i_k} \ran_{1\le j\le d,\,
i_k \in I})$ by adding the first row $(b_{i_0},\ldots\hspace{-1pt},b_{i_d})$.

For $A \in R_1(f)_{(a+1)\beta - \beta_0}$ and 
$B \in R_1(f)_{(b+1)\beta - \beta_0}$ there is a unique constant
$c$ such that 
$$A \cdot B\, (x_1\ldots x_n) - cJ \:\in\: \left\lan x_1 \frac{\partial f}{
\partial x_1},\hspace{2pt}\ldots , x_n \frac{\partial f}{ \partial x_n }\right\ran.$$
Then $$ \int_X {\rm Res}(\omega_A)^{b,a} \cup {\rm Res}(\omega_B)^{a,b}
\;=\; c(-2\pi i)^d \, c_{ab} \, d! \,{\rm Vol}(\Delta_D),$$ where 
$c_{ab}= \frac{(-1)^{a(a+1)/2 +b(b+1)/2 + a^2 +d -1}}{a!b!} $ and
$\Delta_D$ is the polytope associated to the ample divisor $D$. 

 In the mirror quintic situation, $n=d+1=5$, $e_i=v_i$,
 $\beta=\beta_0 =[\sum_{i=1}^5 D_i]$. $\Delta_D$ is the polytope
 with vertices $(1,0,0,0)$, $(0,1,0,0)$, $(0,0,1,0)$, $(0,0,0,1)$
 and $(-1,-1,-1,-1)$. Note that ${\rm Vol}(\Delta_D)= \frac{5}{4!}$. 
 We compute that\/ $(c_I^\beta)= 625$, 
 and 
\begin{align*}
J\:=\; &\psi(x_1x_2x_3x_4)^5 + \psi(x_1x_2x_3x_5)^5 + \psi(x_1x_2x_4x_5)^5 \\[3pt]
  &+ \psi(x_1x_3x_4x_5)^5 + \psi(x_2x_3x_4x_5)^5  + 25(x_1x_2x_3x_4x_5)^4.
\end{align*}

Now there are two cases to consider.
\begin{itemize}
\item[(i)]  $H^{3,0} \cup H^{0,3}\,$:
Both $H^{3,0}(X) \cong R_1(f)_0$ and $H^{1,2}(X) \cong R_1(f)_{3\beta}$ 
 have rank $1$. Choose $A=1 \in R_1(f)_0$ and 
 $B'=(x_1x_2x_3x_4x_5)^3 \in R_1(f)_{3\beta}$. Then we use the Division
 Algorithm in \cite{[CLO]}, p.$\,$61, to compute that\/ $c=\frac{125}{\psi^5 +5^5}\,$. 
 Here we have $a=0$, $b=3$ and thus $c_{ab}=- \frac{1}{3!}$.
 Hence we have $${\rm Res}(\omega_{A})^{3,0} \cup {\rm Res}(\omega_{B})^{0,3} \; =\;
 - \frac{5000 \pi^4}{3(\psi^5 + 5^5)} PD([X]).$$
\item[(ii)] $H^{2,1} \cup H^{1,2}\,$:
Again both $H^{2,1}(X) \cong R_1(f)_\beta$ and $H^{1,2}(X) \cong R_1(f)_{2\beta}$ 
 have rank $1$. Choose $A'=x_1x_2x_3x_4x_5 \in R_1(f)_\beta$ and 
 $B'=(x_1x_2x_3x_4x_5)^2 \in R_1(f)_{2\beta}$. Since $A'\cdot B'$ is same as
 $A\cdot B$, we have the same\/ $c=\frac{125}{\psi^5 +5^5}\,$. But now $a=1$ , 
 $b=2$ and thus  $c_{ab}=\frac{1}{2!}$.
 Hence we have $${\rm Res}(\omega_{A'})^{2,1} \cup {\rm Res}(\omega_{B'})^{1,2} \;=\;
  \frac{5000 \pi^4}{\psi^5 + 5^5} PD([X]).$$
\end{itemize}

\section{$3$-Point Function via Localization }\label{sec:localization}

Consider the general case, i.e. $X$ is a nondegenerate Calabi-Yau hypersurface of a projective
simplicial toric variety $Y$. Then a tricyclic sector $X_{(\bf{g})}$ is a hypersurface of the
corresponding tricyclic sector $Y_{(\bf{g})}$ of $Y$. Denote the associated reduced varieties
by $X_{(\bf{g})}'$ and $Y_{(\bf{g})}'$.
 We showed that $Y_{(\bf{g})}'$ is a toric
variety. There exists a torus equivariant line bundle $L_{(\bf{g})}$ over $Y_{(\bf{g})}'$ with
$X_{(\bf{g})}'$ as the zero locus of a section of $L_{(\bf{g})}$. On the other hand, the 
obstruction
bundle $F_{(\bf{g})}$ over $Y_{(\bf{g})}$ restricts to the reduced obstruction bundle $E_{(\bf{g})}$ over 
$X_{(\bf{g})}$. This follows from the definition of obstruction bundle and the rank formula 
(\ref{eq:dimension-Eg}).
We can associate reduced bundles $E_{(\bf{g})}''$ over $X_{(\bf{g})}'$ to $E_{(\bf{g})}$
and $F_{(\bf{g})}''$ over $Y_{(\bf{g})}'$ to $F_{(\bf{g})}$ by using a $\kG$-invariant homomorphism
of the structure group as in Section 5.  $F_{(\bf{g})}''$ restricts to $E_{(\bf{g})}''$ over $X_{(\bf{g})}'$.
By Poincar\'e duality we have 
\begin{equation}\label{eq:local1}
\int_{X_{(\bf{g})}'} e(E_{(\bf{g})}'') = \int_{Y_{(\bf{g})}'} e(F_{(\bf{g})}'')
\wedge e(L_{(\bf{g})})
=\int_{Y_{(\bf{g})}'} e(F_{(\bf{g})}'' \oplus L_{(\bf{g})}).
\end{equation}

The torus action on $Y_{(\bf{g})}'$ lifts to an action (not unique)
on the bundle $F_{(\bf{g})}'' \oplus L_{(\bf{g})}$ that maps each fiber linearly. Hence we can
 use
the localization technique of Atiyah and Bott \cite{[AB]}, 
adapted to reduced orbifolds, to calculate the 
last integral in (\ref{eq:local1}). 
We carry this out below when $X \subset Y = \cpk /(\zz_5)^3$ 
is a member 
of the mirror quintic family. Let $\Xi$ denote the fan of $\cpk /(\zz_5)^3$ in $N \otimes \rr$ 
as in
Section 2.1 and let $M$ denote the dual lattice of $N$.

Let $g_1 = (\zeta^2 , \zeta^3 , 1, 1, 1)$ and $g_2=g_1$. Then $Y_{(\bf{g})}'$ is isomorphic to 
the toric variety $\{x_1=x_2=0\} \cap Y$. In other words, $Y_{(\bf{g})}'= \bar{O}_{\tau}$ where
$\tau[1] = \{v_1,v_2\}$. Recall that $N_{\tau}$ is the sublattice of $N$ generated by $v_1$ and
$v_2$. Let $N(\tau)=N/N_{\tau}$ be the quotient lattice. The fan of $Y_{(\bf{g})}'$ is given by
the projection of $\Xi$ to $N(\tau)\otimes \rr$. The dual lattice of $N(\tau)$ is $M(\tau) = 
\tau^{\perp} \cap M$.   
The $2$-dimensional torus associated to $Y_{(\bf{g})}'$ is 
${\bf T}= {\rm Spec}(\cc[M(\tau)]) = O_{\tau}$. The characters $\chi^m$ correspond to rational 
functions on $Y_{(\bf{g})}'$ when $m \in M(\tau)$. Let $\{m_1,m_2,m_3,m_4\}$ be the standard
basis of $M$. Then $\{c_1=m_1+m_2+3m_4,\; c_2=m_3-m_4 \}$ is a basis for $M(\tau)$. 
The ${\bf T}$
action on $Y_{(\bf{g})}'$ has three fixed points $q_j=\{x_i = \delta_{i,j}\}, \; 3 \le j \le 5$. 
We study the action of ${\bf T}$ on the normal bundle of $q_j$, i.e., the orbifold tangent space
$(TY_{(\bf{g})}')_{q_j}$. 

Consider $q_5$ first.  Denote local coordinates on a uniformizing sysytem of $Y_{(\bf{g})}'$
around $q_5$ by $z_3=\frac{x_3}{x_5}$ and $z_4 = \frac{x_4}{x_5}$ as in Section 5. Let
$m^1= \frac{1}{5}c_1 + \frac{2}{5}c_2$ and $m^2= \frac{1}{5}c_1 + \frac{1}{5}c_2$.
Then observe that $\lan m^1,v_3 \ran = 1$, $\lan m^1,v_4 \ran = 0$ and $\lan m^1,v_5 \ran =-1$.
Hence $z_3= \chi^{m^1}$ (see Section 3.1). Similarly, $z_4= \chi^{m^2}$. Thus with respect to
the basis $\{c_1,c_2\}$, ${\bf T}$ action is given by
 $(t_1,t_2)(z_3)= t_1^{\frac{1}{5}} t_2^{\frac{2}{5}}z_3$
and  $(t_1,t_2)(z_4)= t_1^{\frac{1}{5}} t_2^{\frac{1}{5}}z_4$. Let $u_1,\, u_2$ be parameters on
 the Lie algebra ${\bf t}_\cc$ of ${\bf T}$ corresponding the above choice of basis. Then the 
${\bf T}$-equivariant Euler class of the normal bundle of $q_5$ is given by  
$e_{\bf{T}}(\nu_{q_5}) =(\frac{1}{5}u_1 + \frac{2}{5}u_2)(\frac{1}{5}u_1 + \frac{1}{5}u_2)$. 
One similarly obtains $e_{\bf{T}}(\nu_{q_4}) 
=(\frac{1}{5}u_2)(-\frac{1}{5}u_1 - \frac{1}{5}u_2)$,
and $e_{\bf{T}}(\nu_{q_3}) =(-\frac{1}{5}u_2)(-\frac{1}{5}u_1 - \frac{2}{5}u_2)$.

Now we want to lift the ${\bf T}$ action to the line bundles $L_{(\bf{g})}$ and $F_{(\bf{g})}''$.
Quite generally, suppose $L$ is a line bundle on a toric variety, corresponding to a Cartier 
divisor $\{U_\sigma, \chi^{-m_\sigma} \}$. The transition functions $h_{\tau \sigma} : U_\sigma 
\times \cc \supset U_{\sigma \cap \tau} \times \cc \to U_{\sigma \cap \tau} \times \cc \subset
U_\tau \times \cc$ for $L$ are given by 
$h_{\tau \sigma}(x,c) = (x, \chi^{(m_\sigma - m_\tau)}(x)c)$. Then one can define a ${\bf T}$
action on $L$ that makes it a ${\bf T}$-equivariant bundle (cf.$\;$\cite{[Od]})
as follows: 
$ t(x,c)= (tx,\chi^{-m_\sigma}(t)c)$ for $t\in {\bf T}$ and $(x,c) \in U_\sigma \times \cc$.

Let $D_i$ denote the Weil divisor $\{ x_i =0\}$ of $Y$. Then the hypersurface $X$ corresponds
to a section of the anticanonical bundle $-K_Y = \sum_{i=1}^5 D_i$. Thus $X\cap \{x_1=x_2=0\}$
corresponds to a section of $-K_Y D_1 D_2$ which is linearly equivalent to $5 D_1 D_2 D_3$. 
Hence $L_{(\bf{g})}=[5 D_1 D_2 D_3]$.
Consider the open covering of $Y_{(\bf{g})}'$ by 
 $U_i \cap Y_{(\bf{g})}' = \{x_1=x_2=0,\, x_i \neq 0 \}, \; 3 \leq i
\leq 5$. $U_i \cap Y_{(\bf{g})}'$
is the affine open set of $Y_{(\bf{g})}'$ corresponding to the cone generated by the projection
 of $\{v_j,v_k\,: \, 3 \le j,k \le 5,\, j,k \neq i \}$. For instance $U_5 \cap Y_{(\bf{g})}'$ is
 generated by projection of $\{v_3,v_4\}$. The equations $\lan -m, v_3 \ran =5$ and 
 $\lan -m, v_4 \ran =0$ have solution $-m = c_1 + 2c_2$ in $M(\tau)$. Hence the divisor 
 $5 D_1 D_2 D_3$ is given by the rational function $\chi^{(c_1+2c_2)}$ on
 $U_5 \cap Y_{(\bf{g})}$.
 Similarly, it is given by the rational functions $\chi^{c_2}$ on $U_4 \cap Y_{(\bf{g})}'$
 and $\chi^0$ on $U_3 \cap Y_{(\bf{g})}'$. Hence the action of ${\bf T}$ on $L_{(\bf{g})}$
 at the fixed points $q_5,\,q_4\;{\rm and}\; q_3$ has weights $(u_1 + 2u_2),\; u_2\;
 {\rm and}\; 0$ respectively. 

The orbifold line bundle $F_{(\bf{g})}$ is trivialized by the generator 
$\frac{\partial}{\partial (x_2/x_i)} \otimes (\zeta \bar{\omega}_1 +\bar{\omega}_2)$ on 
$U_i \cap Y_{(\bf{g})}$.  The transition maps $h_{ij} : (U_j \cap Y_{(\bf{g})}')
\times \cc \supset (U_j \cap U_i \cap Y_{(\bf{g})}') \times \cc \to (U_j \cap U_i
\cap Y_{(\bf{g})}') \times \cc \subset
(U_i \cap Y_{(\bf{g})}') \times \cc\,$ for $F_{(\bf{g})}''$ are given by 
$h_{ji}(x,c) = (x,(\frac{x_j}{x_i})^5(x)c)$. Thus we can define a
 ${\bf T}$ action on $F_{(\bf{g})}''$ by:
\begin{itemize}
\item[(i)] $ t(x,c)= (tx,c) = (tx,\chi^{0}(t)c)$ for $t\in {\bf T}$ and $(x,c) 
\in (U_5\cap Y_{(\bf{g})}') \times \cc$.
\item[(ii)] $ t(x,c)= (tx,(\frac{x_5}{x_4})^5(t)c)=(tx,\chi^{- c_1 - c_2}(t)c)$ for $t\in {\bf T}$ and $(x,c) 
 \in (U_4\cap Y_{(\bf{g})}') \times \cc$.
\item[(iii)] $ t(x,c)= (tx,(\frac{x_5}{x_3})^5(t)c)=(tx,\chi^{- c_1 -2 c_2}(t)c) $ for $t\in {\bf T}$ and $(x,c) 
 \in (U_3\cap Y_{(\bf{g})}') \times \cc$.
\end{itemize}

This action of ${\bf T}$ on $F_{(\bf{g})}''$ at the fixed points 
$q_5,\,q_4\;{\rm and}\; q_3$ has weights $ 0,\;(-u_1 -u_2)\; 
 {\rm and}\; (- u_1 - 2 u_2 )$ respectively. 

Applying corollary 9.1.4 of \cite{[CK]}, 
which shows that for reduced orbifolds the localization formula
has to be modified by dividing the contribution of each fixed 
point by the order of its local group, we have
\begin{align*}
& \int_{Y_{(\bf{g})}'} e(F_{(\bf{g})}'' \oplus L_{(\bf{g})} ) \\[3pt]
 &= \frac{1}{25}
\left(\frac{0(u_1+2u_2)}{(\frac{1}{5}u_1 +\frac{2}{5} u_2)(\frac{1}{5}u_1+\frac{1}{5} u_2) } + 
\frac{(-u_1 - u_2)u_2}{(\frac{1}{5}u_2)(-\frac{1}{5}u_1 -\frac{1}{5} u_2)} +
\frac{(-u_1 - 2 u_2)0}{(-\frac{1}{5}u_2)(-\frac{1}{5}u_1 -\frac{2}{5} u_2)}
\right)\\[3pt]
&=\; 1\; .
\end{align*}

Finally by (\ref{eq:A=A'}) and (\ref{eq:A'=A"}) we must have 
$$\int_{X_{(\bf{g})}}^{orb} e(E_{(\bf{g})})\: = \:\frac{1}{25}\; .$$


\bigskip

\begin{thebibliography}{[CLO]}

\bibitem[AB]{[AB]}  M.F. Atiyah and R. Bott:  The moment map and equivariant 
cohomology, {\it Topology\/} {\bf 23} (1984), 1--28.  

\bibitem[Bai]{[Bai]} W. Baily, Jr.: The decomposition theorem for $V$-manifolds,
{\it Amer. J. Math.\/} {\bf 78} (1956), 862--888. 

\bibitem[Bat]{[B1]} V. Batyrev: Dual polyhedra and mirror symmetry for
 Calabi-Yau hypersurfaces in toric varieties, {\it J. Algebraic Geom.\/}
 {\bf 3} (1994), 493--535; also available at  math.AG/9310003.

\bibitem[BC]{[BC]} V. Batyrev and D. Cox:  On the Hodge structure of
 projective  hypersurfaces in toric varieties, {\it Duke Math. J.\/}
 {\bf 75} (1994), 293--338; also available at math.AG/9306011.

\bibitem[BM]{[BM]} L. Borisov and A. Mavlyutov:  String cohomology of Calabi-Yau hypersurfaces
via Mirror Symmetry, math.AG/0109096.

\bibitem[CG]{[CG]}  J. Carlson and P. Griffiths:  Infinitesimal variations of 
Hodge structure and the global Torelli problem,  
{\sl Jorne\'es de g\'eom\'etrie alg\'ebraic d'Angers}, 1979, 
Sijthoff and Nordhoff, Alphen aan den Rijn, 1980, pp. 51--76. 

\bibitem[CR1]{[CR1]}  W. Chen and Y. Ruan:  A new cohomology theory for orbifold, 
math.AG/0004129 v3.

\bibitem[CR2]{[CR2]}  W. Chen and Y. Ruan:  Orbifold Gromov-Witten theory,  
math.AG/0103156. 

\bibitem[CK]{[CK]}  D.A. Cox and S. Katz:  {\sl Mirror Symmetry and Algebraic
Geometry}.\/  Mathematical Surveys and Monographs, 68.  Amer. Math. Soc., 
Providence, RI, 1999. 

\bibitem[CLO]{[CLO]}  D. Cox, J. Little and D. O'Shea:
{\sl Ideals, varieties, and algorithms.  An introduction to
computational algebraic geometry and commutative algebra.}\/
Second edition. Springer-Verlag, New York, 1997.

\bibitem[Fu]{[Fu]}  W. Fulton:{\sl Introduction to Toric Varieties}.\/
Annals of Mathematics Studies, 131. The William H. Roever Lectures in Geometry.
Princeton Univ. Press, Princeton, NJ, 1993.

\bibitem[GM]{[GM]} M. Goresky and R. MacPherson:
Intersection homology theory,
{\it Topology} {\bf 19} (1980), no. 2, 135--162.

\bibitem[GP]{[GP]} B. Greene and M. Plesser:  Duality in Calabi-Yau
 moduli space, {\it Nuclear Physics} {\bf B338} (1990), 15--37.

\bibitem[Gr]{[Gr]}  P.A. Griffiths:
{\sl Introduction to Algebraic Curves}.\/
Translated from the Chinese by Kuniko Weltin. Translations of Mathematical Monographs, 76.
Amer. Math. Soc., Providence, RI, 1989.

\bibitem[Ka]{[Ka]}  T. Kawasaki:  The signature theorem for $V$-manifolds,
{\it Topology\/} {\bf 17} (1978), 75--83.

\bibitem[KN]{[KN]}  S. Kobayashi and K. Nomizu: 
{\sl Foundations of Differential Geometry, Vol I.}\/ 
Interscience Publishers, a division of John Wiley \& Sons, New York, 1963. 

\bibitem[Ma]{[Ma]}  A.R. Mavlyutov:
Semiample hypersurfaces in toric varieties, {\it Duke Math. J.}\/
{\bf 101} (2000), no. 1, 85--116; also available at
math.AG/9812163 v2.

\bibitem[Od]{[Od]}  T. Oda:  {\sl Convex Bodies and Algebraic Geometry}.\/ 
Springer-Verlag, Berlin-Heidelberg-New York, 1988.

\bibitem[Po]{[Po]}  M. Poddar:  Orbifold Hodge numbers of Calabi-Yau hypersurfaces,
{\it Pacific J. Math.\/} (to appear), available at math.AG/0107152.

\bibitem[Ru]{[Ru2]}  Y. Ruan:  Cohomology ring of crepant resolutions of orbifolds,
math.AG/0108195.

\bibitem[Sc]{[Sc]}  P. Scott: The geometries of $3$-manifolds,
 {\it Bull. London Math. Soc.\/} {\bf 15} (1983), no. 5, 401--487.

\bibitem[Wo]{[Wo]}  S. Wolfram:  {\sl The Mathematica$^{\circledR}$ Book}.\/
Fourth edition. Wolfram Media, Inc., Champaign, IL; Cambridge
University Press, Cambridge, 1999.

\end{thebibliography}
\end{document}